\documentclass[12pt]{article}
\usepackage{amsmath}
\usepackage{amssymb}
\usepackage{amsthm}
\usepackage{amscd}
\usepackage{amsxtra}
\usepackage{verbatim}
\usepackage{xcolor}
\usepackage{color}
\usepackage{enumerate}
\usepackage{mathrsfs}
\usepackage[all]{xy} 

\usepackage{array,arydshln}
\usepackage{bm}

\allowdisplaybreaks

\numberwithin{equation}{section}

\definecolor{dblue}{rgb}{0,0,0.45}
\definecolor{red}{rgb}{0.7,0,0}

 \RequirePackage{geometry}
\newtheorem{theorem}{Theorem}[section]
\newtheorem{lemma}[theorem]{Lemma}
\newtheorem*{lemma*}{Lemma}
\newtheorem{corollary}[theorem]{Corollary}
\newtheorem{proposition}[theorem]{Proposition}

\theoremstyle{definition}

\newtheorem{remark}[theorem]{Remark}

\newtheorem{example}[theorem]{Example}

\theoremstyle{remark}

\newcommand{\N}{{\mathbb N}}
\newcommand{\R}{{\mathbb R}}

\newcommand{\E}{{\mathbb E}}
\newcommand{\cA}{{\mathcal A}}
\newcommand{\cC}{{\mathcal C}}

\newcommand{\cH}{{\mathcal H}}

\newcommand{\cK}{{\mathcal K}}
\newcommand{\cL}{{\mathcal L}}

\newcommand{\cV}{{\mathcal V}}
\newcommand{\cW}{{\mathcal W}}

\newcommand{\al}{\alpha}

\newcommand{\la}{\langle}
\newcommand{\ra}{\rangle}
\newcommand{\nn}{\nonumber}

\newcommand{\vertiii}[1]{{\left\vert\kern-0.25ex\left\vert\kern-0.25ex\left\vert #1 
    \right\vert\kern-0.25ex\right\vert\kern-0.25ex\right\vert}}

\date{}
\begin{document}

\title{
Support theorem for pinned diffusion processes
}
\author{  Yuzuru \textsc{Inahama} 
}
\maketitle

\begin{abstract}
\noindent
In this paper we prove a support theorem of 
Stroock-Varadhan type for pinned diffusion processes.
To this end we  use two powerful results from stochastic analysis.
One is quasi-sure analysis for Brownian rough path.
The other is Aida-Kusuoka-Stroock's positivity theorem 
for the densities of weighted laws of non-degenerate 
Wiener functionals.
\end{abstract}

\section{Introduction}

Let us consider the following Stratonovich stochastic 
differential equation (SDE) on $\R^e$ ($e \ge 1$)
driven by a standard $d$-dimensional
 Brownian motion $w=(w_t)_{0\le t \le 1}$:
\[
dX_t = \sum_{i=1}^d  V_i ( X_t)\circ dw_t^i + V_0 ( X_t) dt,
\qquad
 X_0 =a \in \R^e.
 \]
Here, $V_i~(0\le i \le d)$ are sufficiently nice vector fields on $\R^e$ and
$a \in \R^e$ is an arbitrary (deterministic) starting point.
Throughout this paper the time interval is $[0,1]$ unless otherwise stated.
The corresponding skeleton ordinary differential equation (ODE) 
is given as follows:
For a $d$-dimensional Cameron-Martin path $h\colon [0,1]\to \R^d$, 
\[
dx^h_t = \sum_{i=1}^d  V_i ( x^h_t) dh_t^i + V_0 ( x^h_t) dt,
\qquad
 x^h_0 =a \in \R^e.
\]
We will write $\Psi (h) = (x^h_t)_{0\le t \le 1}$ for simplicity
and denote by ${\cal H}$ the set of all $d$-dimensional
Cameron-Martin paths.

We are interested in the (topological) support of the law of 
the diffusion process $X=(X_t)_{0\le t \le 1}$
and would like to describe it in terms of the skeleton ODE.
Stroock-Varadhan's support theorem states that 
the support equals the closure of $\{ \Psi (h) \colon h \in {\cal H}\}$. 
(See \cite{sv} and \cite[Section 8.3]{strbk}.)
In the original work the uniform topology was used,
but later it was improved to the $\alpha$-H\"older topology 
with $0 <\alpha<1/2$ in \cite{ms,bagl}.
A quite general approach to support theorems
of this kind  in \cite{aks} should also be referred to.


After the pioneering work \cite{sv}, the support theorem 
became one of central topics in the study of SDEs
and were generalized to many directions.
A partial list could be as follows.
A generalization to SDEs with unbounded coefficients was done in \cite{gp}.
Support theorems for reflecting diffusions were proved in \cite{dp, rw}.
The topology of the path space was further refined in \cite{gnss}.
The case of anticipating SDEs were studied in \cite{mn_anti, ms_anti}.
The case of (Volterra-type) SDEs with path-dependent coefficients 
were recently studied in \cite{ck, kal}.
A support theorem for McKean-Vlasov SDEs was proved in \cite{xg}.
A support theorem for jump-type SDEs was studied in \cite{simon}.
For support theorems for stochastic PDEs, see \cite{ms_spde1, bms, ms_spde2, nak} among others.
(Results related to rough path theory will be listed shortly.)


Using rough path theory, 
Ledoux-Qian-Zhang \cite{lqz} gave a new proof to the support theorem
twenty years ago.
Their idea could be summarized as follows.
If the It\^o map $w \mapsto X$, i.e. the solution map of 
the above SDE, were continuous, then the proof of the support 
theorem would be simple. 
(In fact, it is not continuous.  So, the proof is not easy.)
Compared to the usual SDE theory,  rough path theory
has a prominent feature.
The Lyons-It\^o map $\Phi$, i.e. the solution map of the corresponding 
rough differential equation (RDE) is continuous.
Moreover,  $X = \Phi ({\bf W})$, a.s. and $\Phi$ is compatible with $\Psi$.
Here, ${\bf W}$ is Brownian rough path, i.e. the standard 
Stratonovich rough path lift of $w$.
Hence, if a support theorem for the law of ${\bf W}$ is obtained
on the geometric rough path space, 
Stroock-Varadhan's support theorem follows immediately.
The support theorem for ${\bf W}$ was first proved 
with respect to the $p$-variation topology ($2<p<3$) in \cite{lqz}
and then improved to the case of the
$\alpha$-H\"older topology ($1/3 <\alpha <1/2$) in \cite{friz}.
This support theorem was later generalized for the laws of 
Gaussian rough paths in \cite{fv10}.  
(See also \cite[Sections 13.7 and 15.8]{fvbk} and reference therein.)
Other applications of rough path technique to support theorems 
are found in \cite{dd, aida_pre, crs}.
Support theorems are also studied  
in the theory of singular stochastic PDEs, 
which is a descendant of rough path theory. See \cite{cf, tw, hs, mats}.


In this paper we study an analogous support theorem 
for the law of the pinned diffusion process 
which is condition to end at $b \in \R^e$ at the time $t=1$, i.e. 
the law of $X=(X_t)_{0\le t \le 1}$ under the conditional measure 
${\mathbb E}[ \,\cdot\, \mid X_1 =b\,]$ (heuristically).
A natural guess could be that  its support equals the closure of 
$\{ \Psi (h) \colon h \in {\cal H}, \Psi (h)_1=b\}$. But, is it really true?

Before discussing this problem, 
we first review a positivity theorem \cite{bal, aks}
for the density of the law of $X_t$, which is closely related to 
the support theorem.
Under a H\"ormander-type condition on $V_i$'s, 
the law of $X_t$ has a smooth density $p(t,a,y)$
with respect to the Lebesgue measure for every $t \in (0,1]$.
It is natural and important to ask whether  
or under what condition $p(t,a,y)>0$.
(For instance, if $p(1,a,b)=0$, the above-mentioned pinned diffusion 
measure does not exist.)
The positivity theorem states that $p(t,a,b)>0$
if and only if there exists $h \in {\cal H}$ such that
$\Psi (h)_t=b$ and $D\Psi (h)_t \colon {\cal H} \to \R^e$ 
is a surjective linear map.
Here, $D$ stands for the Fr\'echet derivative on ${\cal H}$.
The first paper which proved this result was \cite{bal}. 
Then, a very general result by Aida-Kusuoka-Stroock \cite{aks} followed,
which will be used in this paper.

A significant feature of \cite{aks} is that it studies the positivity 
of the density of a {\it weighted} law of a Wiener functional.
(In most of the works on this problem, the weight identically 
equals $1$.)
In the proof of our main theorem, we will exploit this
arbitrariness of the weight.
To be more specific, we will choose as a weight
a Wiener functional that looks like the indicator function 
of an open neighborhood of a given geometric rough path.


If we keep these two famous theorems in mind, 
we can guess what the support of the pinned 
diffusion measure looks like.
First, let us first recall a precise definition 
of the pinned diffusion measure $\mathbb{Q}_{a,b}$ ($a, b \in \R^e$).
We assume that $V_i ~(0 \le i \le d)$ satisfy 
H\"ormander's bracket generating condition at every $x \in \R^e$
(see Remark \ref{rem.nondeg} (A)).
Then, the density $p (t,x,y)$ exists for all $x,y \in \R^e$ and $t\in(0,1]$.
We further assume that $p (1,a,b) >0$, which is equivalent to 
the existence of $h \in {\cal H}$ such that
$\Psi (h)_1=b$ and $D\Psi (h)_1 \colon {\cal H} \to \R^e$ 
is surjective.
For every $\beta \in (1/3, 1/2)$, 
$\mathbb{Q}_{a,b}$ is a unique probability measure on 
the $\beta$-H\"older continuous path space
\[
\cC_{a,b}^{\beta\textrm{-H}} (\R^e)
:= \{\xi \colon [0,1]\to \R^e 
\colon \mbox{$\beta$-H\"older continuous and }
\xi_0 =a, \, \xi_1 =b \}
\]
 with the following property (it does exist):
For every $k\ge 1$,  $\{0 <t_1<\cdots < t_k <1\}$
and $g_1, \ldots, g_k \in C_0^\infty (\cV)$, 
\begin{eqnarray*}
\int   \prod_{i=1}^k g_i (\xi_{t_i})
\mathbb{Q}_{a,b} (d\xi)
&=&
p (1, a,b)^{-1} 
\int_{(\R^e)^k}  \Bigl\{ \prod_{i=1}^k g_i (z_i) \Bigr\} p (t_1, a,z_1)  
\nn\\
&& \times
\Bigl\{\prod_{i=2}^{k}p (t_i - t_{i-1}, z_{i-1}, z_{i})  \Bigr\}
 p (1- t_k, z_{k},b) 
\Bigl\{\prod_{i=1}^k dz_i \Bigr\}.
\end{eqnarray*}
Here, $dz_i$ ($1 \le i \le k$) is the Lebesgue measure on $\R^e$ and 
$(\xi_t)$ is the canonical coordinate process 
on $\cC_{a,b}^{\beta\textrm{-H}} (\R^e)$.

We will prove in Corollary \ref{cor.supp_X} that 
the support of $\mathbb{Q}_{a,b}$ equals the closure 
with respect to the $\beta$-H\"older topology of 
\[
\{ \Psi (h)  \colon   \mbox{$h \in \cH$,  \, $\Psi (h)_1=b$, \,
$D \Psi (h)_1
\colon \cH\rightarrow \mathbb{R}^e$ is surjective} 
\}.
\]
In fact, this is a special case  of our more general result
(Corollary \ref{cor.supp_Y}), in which we will prove
a support theorem for generalized pinned measure.
With quasi-sure analysis and Malliavin calculus,
 one can easily see that this kind of generalized pinned measures exist.
However, since it is difficult to give a brief introduction of them, 
we do not explain Corollary \ref{cor.supp_Y} here.

These two corollaries are direct consequences of our 
main theorem (Theorem \ref{or.supp.GRP}).
By a well-known theorem in quasi-sure analysis, 
there exists a measure on the classical Wiener space that looks 
like a pullback of $\mathbb{Q}_{a,b}$ by the It\^o map.
Since the rough path lift map is in fact quasi-surely defined, 
the measure admits a lift to a measure on the
geometric rough path space. 
By $\infty$-quasi-continuity of the lift map,
its image measure
 induced by the Lyons-It\^o map is $\mathbb{Q}_{a,b}$ as expected.
Theorem \ref{or.supp.GRP} is a support theorem 
for this lifted measure.
To show it, we use quasi-sure analysis and 
Aida-Kusuoka-Stroock's positivity theorem \cite[Theorem 2.8]{aks}.
For precise formulations and statements of
 these  results, see Section \ref{sec.SUP}.


The organization of this paper is as follows.
Section \ref{sec.prelim} is devoted to reviewing known results 
of Malliavin calculus that will be used in the main part of this paper.
After we recall fundamentals of (Watanabe's distributional) 
Malliavin calculus and quasi-sure analysis, 
we review Aida-Kusuoka-Stroock's positivity theorem, which will play a
major role in our proof.
In Section \ref{sec.RP} we recall basic facts on
quasi-sure analytic properties of Brownian rough path.
In relation to this, a Besov-type topology is introduced on the geometric 
rough path space. 
Section \ref{sec.cK_diff} is a core part of this work, in which
we prove twice $\cK$-differentiability of the Lyons-It\^o map,
that is, the solution map of a rough differential equation.
This property is the key condition in Aida-Kusuoka-Stroock's theorem.
In Section \ref{sec.SUP} we state our main theorems precisely 
and prove them rigorously.
Our key result is Theorem \ref{or.supp.GRP}.
This is a support theorem on a geometric rough path space 
for a measure that looks like the ``pullback" by the Lyons-It\^o map
of a pinned diffusion measure.
Since the Lyons-It\^o map is continuous, the
support theorems for pinned diffusion measures
(Corollaries \ref{cor.supp_Y} and \ref{cor.supp_X}) follow immediately.

\medskip

\noindent 
{\bf Notation:} In the sequel we will use the following notation.
We write $\N =\{1,2, \ldots\}$.
The time interval of (rough) paths 
and stochastic processes is $[0,1]$ throughout the paper.
Below, we assume $d \in \N$.

\begin{itemize} 
\item
The set of all continuous paths $\varphi\colon [0,1] \to\R^d$ 
is denoted by $\cC (\R^d)$. 
With the usual sup-norm $\|\varphi\|_{\infty} := \sup_{0\le t \le 1}|\varphi_t|$ 
on the $[0,1]$-interval,
this is a Banach space.
The increment of $\varphi$ is often denoted by $\varphi^1$,
that is, $\varphi^1_{s,t} := \varphi_t - \varphi_s$ for $s\le t$.
For $a, b \in \R^d$, we write
$\cC_a (\R^d) =\{ \varphi \in \cC (\R^d) \colon \varphi_0 =a\}$
and, in a similar way, $\cC_{a,b} (\R^d) =\{ \varphi \in \cC (\R^d) \colon \varphi_0 =a, \varphi_1 =b\}$.

\item
Let $\alpha \in (0,1]$. The $\alpha$-H\"older seminorm of 
$ \varphi \in \cC (\R^d)$ is defined as usual by
\[
\|  \varphi\|_{\alpha} :=\sup_{0\le s<t \le 1} \frac{|\varphi^1_{s,t}|}{(t-s)^{\alpha}}.
\]
The $\alpha$-H\"older continuous path space is 
denoted by $\cC^{\alpha\textrm{-H}} (\R^d) =\{ \varphi \in \cC (\R^d)
\colon  \|  \varphi\|_{\alpha} <\infty\}$, which is a non-separable
Banach space
with the norm $\|  \varphi\|_{\alpha}  +|\varphi_0|$.
The closure of 
$\{\varphi \in \cC (\R^d) \colon \mbox{$\varphi$ is piecewise-$C^1$} \}$
with respect to the $\alpha$-H\"older topology is denoted by 
$\cC^{0, \alpha\textrm{-H}} (\R^d)$.
This is a separable Banach subspace of $\cC^{\alpha\textrm{-H}} (\R^d)$.
For a starting point $a\in \R^d$ and 
an end  point $b\in \R^d$, $\cC_a^{0, \alpha\textrm{-H}} (\R^d)$
and $\cC_{a,b}^{0, \alpha\textrm{-H}} (\R^d)$
are defined in an analogous way as above.

\item
For $1/3 <\alpha \le 1/2$, 
$G\Omega^{\textrm{H}}_{\alpha} ( {\mathbb R}^d)$
stands  for  the $\alpha$-H\"older 
geometric rough path space over ${\mathbb R}^d$.
A generic element of $G\Omega^{\textrm{H}}_{\alpha} ( {\mathbb R}^d)$
is denoted by ${\bf w} =({\bf w}^1, {\bf w}^2)$.
(See \cite{fvbk,lcl} among others for basic information on 
geometric rough paths.)

\item
The Cameron-Martin space associated with standard  $d$-dimensional
Brownian motion is denoted by $\cH =\cH^d$
(except in Section \ref{sec.prelim}).
Its precise definition is given as follows:
$
\cH := \{  h \in \cC_0 (\R^d) \colon \mbox{$h$ is absolutely continuous
and $\|h \|_{\cH} <\infty$} \}$, 
where 
\[
\langle h,k\rangle_{\cH}:= 
   \int_0^1 \langle h_s^\prime, k_s^\prime \rangle_{\R^d} ds 
\quad \mbox{and}\quad 
\|h \|_{\cH} := \langle h,h\rangle_{\cH}^{1/2},
\qquad h, k \in \cH.
\]
This is a real separable Hilbert space with this inner product.
It is easy to see that $\cH \subset \cC_0^{0, 1/2\textrm{-H}} (\R^d)$.
If we set ${\cal L} (h)^1_{s,t} =h_t-h_s$ and 
${\cal L} (h)^2_{s,t} =\int_s^t   (h_u-h_s) \otimes dh_u$ 
for $h\in \cH$ and $0\le s \le t \le 1$, 
then ${\cal L}\colon \cH \hookrightarrow G\Omega^{\textrm{H}}_{1/2} ( {\mathbb R}^d)$
becomes a locally Lipschitz continuous injection.
(A map between two metric space is said to be 
locally Lipschitz continuous if the map, when restricted to every bounded subset of the domain, is Lipschitz continuous.)  
We call ${\cal L} (h)$ the natural lift of $h$ and 
will sometimes denote it by ${\bf h}$.

\item
Let $U$ be an open subset of $\R^m$.
For $k \in \N \cup \{0\}$,  $C^k (U, \R^n)$ denotes the set of 
$C^k$-functions from $U$ to $\R^n$.
(When $k=0$, we simply write $C (U, \R^n)$ 
instead of $C^0 (U, \R^n)$.)
The set of bounded $C^k$-functions $f \colon U\to \R^n$
whose derivatives up to order $k$ are all bounded 
is denoted by $C_b^k (U, \R^n)$. This is a Banach space with 
$\| f\|_{C_b^k } := \sum_{i=0}^k \|\nabla^i f\|_{\infty}$.
(Here, $ \|\cdot\|_{\infty}$ stands for the usual sup-norm on $U$.)
As usual, we set $C^\infty (U, \R^n) := \cap_{k =0}^\infty C^k (U, \R^n)$
and $C^\infty_b (U, \R^n) := \cap_{k=0}^\infty C^k_b (U, \R^n)$.
\end{itemize}

%
\section{Preliminaries from Malliavin calculus}
\label{sec.prelim}

In this section,  $(\cW, \cH, \mu)$ is 
an abstract Wiener space. 
That is, ~$(\cW, \|\cdot\|_{\cW})$~is a separable Banach space, ~$(\cH, \|\cdot\|_{\cH})$~ is a separable Hilbert space, $\cH$~is a dense subspace of~$\cW$~and the inclusion map is continuous, and~$\mu$~is a (necessarily unique) Borel probability measure
 on $\cW$
 with the property that
\begin{equation}\label{abspro}
\int_{\cW}\exp\Bigl(\sqrt{-1}_{\cW^*}\langle \lambda, w \rangle_{\cW}\Bigr)\mu (d w)=\exp\Bigl(-\frac{1}{2}\|\lambda\|^2_{\cH}\Bigr),
\qquad
\qquad
\lambda \in \cW^* \subset \cH^*,
\end{equation}
where we have used the fact that $\cW^*$ becomes a dense subspace of $\cH$ when we make the natural identification between $\cH^*$ and $\cH$ itself. 
Hence, $\cW^* \hookrightarrow \cH^* =\cH \hookrightarrow\cW$
and both inclusions are continuous and dense.
We denote by $\{ \langle k, \bullet \rangle \colon k\in \cH\}$ 
the family of centered Gaussian random variables 
defined on $\cW$ indexed by $\cH$
(i.e. the homogeneous Wiener chaos of order $1$).
If $\langle k, \bullet \rangle_{\cH} \in \cH^*$ 
extends to an element of $\cW^*$,
then the extension coincides with the  random 
variable $\langle k, \bullet \rangle$.
(When $\langle k, \bullet \rangle_{\cH} \in \cH^*$ does not 
extend to an element of $\cW^*$, $\langle k, \bullet \rangle$
is define as  the $L^2$-limit of $\{\langle k_n, \bullet \rangle \}_{n=1}^\infty$, where  $\{k_n\}_{n=1}^\infty$ 
is any sequence of $\cH$ such that 
$\langle k_n, \bullet \rangle_{\cH}\in \cW^*$ for all $n$ and $\lim_{n\to\infty} \|k_n -k\|_\cH=0$.)
We also denote by $\tau_k \colon \cW \to \cW$ the translation 
$\tau_k (w) =w+k$.
(For basic information on abstract Wiener spaces, see \cite{sh, hu} among others.)

\subsection{Watanabe distribution theory
and quasi-sure analysis}
We first quickly summarize some basic facts in Malliavin calculus,
which are related to Watanabe distributions
(i.e. generalized Wiener functionals) 
and quasi-sure analysis.
Most of the contents and the notation
in this section are found in 
 \cite[Sections V.8--V.10]{iwbk} with trivial modifications.
Also, \cite{sh, nu, hu, mt, kunita} are good textbooks of Malliavin calculus.
For basic results of quasi-sure analysis, we refer to \cite[Chapter II]{ma}.
We use the following notation and facts in the main part of this paper.
\begin{enumerate}
\item[{\bf (a)}]
Sobolev spaces ${\bf D}_{p,r} ({\cal K})$ of ${\cal K}$-valued 
(generalized) Wiener functionals, 
where ${\cal K}$ is a real separable Hilbert space and
$p \in (1, \infty)$, $r \in {\mathbb R}$.
As usual, we will use the spaces 
${\bf D}_{\infty} ({\cal K})= \cap_{k=1 }^{\infty} \cap_{1<p<\infty} {\bf D}_{p,k} ({\cal K})$, 
$\tilde{{\bf D}}_{\infty} ({\cal K}) 
= \cap_{k=1 }^{\infty} \cup_{1<p<\infty}  {\bf D}_{p,k} ({\cal K})$ of test functionals 
and  the spaces ${\bf D}_{-\infty} ({\cal K}) = \cup_{k=1 }^{\infty} \cup_{1<p<\infty} {\bf D}_{p,-k} ({\cal K})$, 
$\tilde{{\bf D}}_{-\infty} ({\cal K}) = \cup_{k=1 }^{\infty} \cap_{1<p<\infty} {\bf D}_{p,-k} ({\cal K})$ of 
 Watanabe distributions as in \cite{iwbk}.
When ${\cal K} ={\mathbb R}$, we simply write ${\bf D}_{p, r}$, etc.
\item[{\bf (b)}] 
For 
$F =(F^1, \ldots, F^e) \in {\bf D}_{\infty} ({\mathbb R}^e)$, we denote by 
$\sigma^{ij}_F (w) =  \la DF^i (w),DF^j (w)\ra_{{\cal H}}$
 the $(i,j)$-component of Malliavin covariance 
 matrix ($e\in \N$, $1 \le i,j \le e$).
We say that $F$ is non-degenerate in the sense of Malliavin
if $(\det \sigma_F)^{-1} \in  \cap_{1<p< \infty} L^p (\mu)$.
Here, $D$ is the $\cH$-derivative
(the gradient operator in the sense of Malliavin calculus).
If $F \in {\bf D}_{\infty} ({\mathbb R}^e)$
is non-degenerate, its law on $\R^e$ admits a smooth,
rapidly decreasing
density $p_F =p_F (y)$ with respect to the Lebesgue measure $dy$,
that is, $\mu \circ F^{-1}= p_F (y)dy$.
(This fact is quite famous. See any textbook of Malliavin calculus.)
\item
[{\bf (c)}] Pullback $T \circ F =T(F)\in \tilde{\bf D}_{-\infty}$ of a tempered Schwartz distribution $T \in {\cal S}^{\prime}({\mathbb R}^e)$
on ${\mathbb R}^e$
by a non-degenerate Wiener functional $F \in {\bf D}_{\infty} ({\mathbb R}^e)$. 
The most important example of $T$ is Dirac's delta function.
In that case, ${\mathbb E}[\delta_y (F)] :=\la \delta_y (F),1\ra =p_F (y)$ 
holds for every $y\in\R^e$.
Here, $\la \star, *\ra$ denotes the pairing of ${\bf D}_{-\infty}$ and ${\bf D}_{\infty}$ as usual.
(See \cite[Section 5.9]{iwbk}.)

\item[{\bf (d)}]
This is a continuation of {\bf (b)} and {\bf (c)} above. 
Assume in addition that $G \in {\bf D}_{\infty}$ is non-negative.
Then, $(Gd\mu) \circ F^{-1}$ is called the law of $F$ {\it weighted by} $G$.
(In other words, this law is a probability measure on $\R^e$
determined by $A \mapsto {\mathbb E}[{\bf 1}_A (F) G]$,
where $A$ is a Borel measurable subset of $\R^e$.)
If $F$ is non-degenerate, this law admits a smooth,
rapidly decreasing
density $p_{F,G} =p_{F, G} (y)$ with respect to the Lebesgue measure $dy$,
that is, $(Gd\mu) \circ F^{-1}= p_{F, G} (y)dy$.
In the language of Watanabe distributions, 
we have ${\mathbb E}[\delta_y (F) G]:=\la \delta_y (F),G\ra
 =p_{F,G} (y)$ for every $y\in\R^e$.
(For weighted laws of non-degenerate Wiener functionals, 
we refer to \cite[Sections 5.3 and 5.12]{kunita}.)

\item[{\bf (e)}]
If $\eta \in {\bf D}_{-\infty}$ satisfies that
${\mathbb E}[\eta \, G]:=\la \eta, G \ra \ge 0$ for every non-negative $G \in {\bf D}_{\infty}$,
it is called a positive Watanabe distribution.
According to Sugita's theorem (see \cite{sugita}
or \cite[p. 101]{ma}), for every positive Watanabe distribution $\eta$,
there uniquely exists a finite Borel measure $\mu_{\eta}$ on $\cW$
such that
\[
\la \eta, G\ra = \int_{\cW} \tilde{G} (w)  \mu_{\eta} (dw), 
\qquad G \in {\bf D}_{\infty}
\]
holds,
where $\tilde{G}$ stands for an $\infty$-quasi-continuous modification of 
$G$.
If $\eta \in {\bf D}_{p, -k}$ is positive, then it holds that 
\[
\mu_{\eta}(A) \le \| \eta \|_{p,-k} {\rm Cap}_{q,k} (A)
\qquad
\mbox{for every Borel subset $A\subset \cW$,}
\]
where $p, q  \in (1, \infty)$ with $1/p +1/q =1$, $k \in \N$,
and ${\rm Cap}_{q,k}$ stands for the $(q,k)$-capacity
associated with ${\bf D}_{q,k}$.
(For more details, see \cite[Chapter II]{ma}.)
\end{enumerate}

\begin{remark} 
In some of the books cited in this subsection (in particular \cite{iwbk, ma, mt}),  results are formulated on  a special Gaussian space.
However, almost most all of them 
(at least, those that will be used in this paper) still hold true
 on any abstract Wiener space.
\end{remark}

\subsection{$\cK$-regularity and $\cK$-differentiability}
\label{sec.aks}

In this subsection we quickly review Aida-Kusuoka-Stroock's 
result on the positivity of the density
for non-degenerate Wiener functionals (see \cite{aks}).

Let we first recall the definitions of $\cK$-continuity,
$\cK$-regularity, uniformly $\cK$-regularity, 
and $l$-times $\cK$-regular differentiability, 
which were first introduced in \cite{aks}. 
Note that in these definitions, functions and maps on $\cW$
are viewed as 
everywhere-defined ones (not equivalence classes 
with respect to $\mu$).
It should be noted that
these definitions depend on the choice of exhaustion $\cK$.

For a finite dimensional subspace $K$ of $\cH$,   
$P_K\colon \cH \to K$ stands for the orthogonal projection
and we write $P_K^\perp = {\rm Id}_\cH - P_K$. 
This projection naturally extends to $\bar{P}_K\colon \cW \to K$
as follows:
\[
\bar{P}_K (w) = \sum_{i=1}^{\dim K} \langle e_i, w \rangle e_i,
\]
where $\{ e_i\}_{i=1}^{\dim K}$ is an orthonormal basis of $K$.
(This right hand side
 is independent of the choice of $\{ e_i\}$.)
We set $\bar{P}_K^\perp =  {\rm Id}_\cW - \bar{P}_K$.

Assume that $\cK=\{K_n\}^{\infty}_{n=1}$ is a non-decreasing, countable exhaustion of $\cH$ by finite dimensional subspaces, 
that is, $K_n \subset K_{n+1}$ for all $n$
and $\cup_{n=1}^\infty K_n$ is dense in $\cH$.
Set $P_n=P_{K_n}$, define $\bar P_n$, $P_n^{\perp}$, $\bar P_n^{\perp}$ 
accordingly. We say that a map $F$ from $\cW$ into 
a Polish space $(E, \rho_E)$ is {\it $\cK$-continuous} 
if it is measurable and, 
for each $n\in \N$, 
there is a measurable map $F_n\colon\cW\times K_n\longmapsto E$ with the properties that $F\circ \tau_k=F_n(\cdot,k)$ ($\mu$-a.s.)
 for each $k\in K_n$ and $k\in K_n \longmapsto F_n(w,k)\in E$ is continuous for each $w \in \cW$. Given a $\cK$-continuous 
 map $F$,  we set
\begin{eqnarray}
 F_n^{\perp}(w,k) = F_n(w,-\bar P_n (w) +k)
\qquad {\rm for}~n\in \N ~{\rm and}~k\in K_n.
\end{eqnarray}

Given a measurable map $F \colon\cW\rightarrow E$, we say that $F$ is {\it $\cK$-regular} if $F$ is $\cK$-continuous and there is a continuous map
$\tilde {F} \colon\cH \rightarrow E$ such that
\begin{eqnarray}\label{rhoe}
 \lim_{n \rightarrow \infty} \mu \bigg(\bigg\{w\colon\rho_E(\tilde{F}\circ \bar{P}_n(w),F (w))\vee
  \rho_E(\tilde{F}(h),F_n^{\perp}(w, P_n (h)))
 \geq \epsilon\bigg\}\bigg)=0
 \end{eqnarray}
holds for every $\epsilon>0$ and $h\in \cH$. 
In this case $\tilde {F}$ is called a $\cK$-regularization of $F$.

If  $F$ is a map from $\cW$ into a Polish space $E$, we say that it is {\it uniformly $\cK$-regular} if it is $\cK$-regular and (\ref{rhoe}) can be replaced by the condition that
\begin{eqnarray}
 \lim_{n \rightarrow \infty} \mu \bigg(\bigg\{w &\colon& \sup_{k\in K_m, \|k\|_{\cH}\leq r}\rho_E(\tilde{F}(\bar{P}_n(w)+k),F_{n}(w,k)) \cr
 &&\vee \rho_E(\tilde{F}(h+k), 
 F^{\perp}_{n}(w, P_n(h)+k)) \geq \epsilon\bigg\}\bigg)=0
 \label{rhoe2}
  \end{eqnarray}
for every $m\in \N, r>0, \epsilon>0$ and $h\in \cH$. 
(In \eqref{rhoe2} above and \eqref{rhoe3} below,
we implicitly assume $n \ge m$ since we let $n\to\infty$ for each fixed $m$.)

Let $E$ be a separable Banach space 
and $F$ be a map from $\cW$ into $E$.
Given $l\in \N$ we say that $F$ is 
{\it $l$-times $\cK$-regularly differentiable} 
if $F$ is uniformly $\cK$-regular, $F_n (w,\cdot)$ is $l$-times continuously Fr\'{e}chet differentiable on $K_n$ for each $n\in \N$ and $w \in \cW$, $\tilde{F}$ is $l$-times continuously 
Fr\'{e}chet differentiable 
on $\cH$, and  (\ref{rhoe2}) can be replaced by the condition that
\begin{eqnarray}\label{rhoe3}
 \lim_{n \rightarrow \infty} \mu \bigg(\bigg\{w &\colon& \|\tilde{F}(\bar{P}_n(w)+\bullet)-F_{n}(w,\bullet))\|_{C_{b}^l (B_{K_m}(0,r), E)} \cr
 &&\vee \|\tilde{F}(h+\bullet)-F^{\perp}_{n}(w, P_n(h)+\bullet)\|_{C_{b}^l (B_{K_m}(0,r), E)} \geq \epsilon\bigg\}\bigg)=0
 \end{eqnarray}
for every $m\in \N, r>0, \epsilon >0$ and $h\in \cH$. 
Here, $B_{K_m}(0,r)=\{ k \in K_m \colon \|k\|_{\cH}< r \}$.


The following theorem is \cite[Theorem 2.8]{aks}
(translated into the language of Watanabe distribution theory),
which is the key tool in this paper.
It is a quite general result on the positivity of the density 
function of the law of a non-degenerate Wiener functional.
At first sight, it may not be clear why 
the case of non-constant weight $G$ is so important. 
However, in the proof of our main theorem, the weight will play a crucial role.
%

%
%
\begin{theorem}\label{aks>0}
Let $F\in \mathbf{D}_{\infty} (\mathbb{R}^e)$, $e \in \N$, and
$G\in \mathbf{D}_{\infty}$.
Suppose that $F$ is non-degenerate in the sense of Malliavin
and $G$ is non-negative.
Suppose further that $F$ is twice $\cK$-regularly differentiable 
and $G$ is $\cK$-regular with their $\cK$-regularizations 
$\tilde{F}$ and $\tilde{G}$, respectively. 
Then,  for $y\in \mathbb{R}^e$, the following  are equivalent:
\begin{itemize}
\item
${\mathbb E}[\delta_y (F) G]>0$.
\item
There exists $h \in \cH$ such that 
$D\tilde{F}(h)\colon \cH\rightarrow \mathbb{R}^e$ is surjective, 
$\tilde{F}(h)=y$ and 
$\tilde{G}(h)>0$. 
\end{itemize}
\end{theorem}

\begin{remark}
As is well-known, the condition that ``$D\tilde{F}(h)\colon \cH\rightarrow \mathbb{R}^e$ is surjective" in the above theorem is equivalent to 
non-degeneracy of deterministic Malliavin covariance matrix
of $\tilde{F}$ at $h$.
\end{remark}


The most typical example of $\cK=\{K_n\}$ and $P_n =P_{K_n}$ is 
the dyadic piecewise linear approximation $w(n)$ of the standard 
$d$-dimensional Brownian motion 
$w= (w_t)_{0\le t \le 1}$. 
As usual, $w(n)$ is defined as follows:
$w(n)_{j2^{-n}} =w_{j2^{-n}}$ for all $0\le j \le 2^n$
and $w (n)$
is linearly interpolated on each subinterval 
$[(j-1)2^{-n}, j2^{-n}]$, $0\le j \le 2^n$.

\begin{example}\label{exm.DPL}
Let 
$(\cW, \cH, \mu)$ be the $d$-dimensional 
classical Wiener space, that is, 
{\rm (i)} $\cW := \cC_0 (\mathbb{R}^d)$ is the Banach space of 
$\mathbb{R}^d$-valued
continuous paths that start at $0$ equipped with the usual sup-norm,
{\rm (ii)} $\mu$ is the $d$-dimensional Wiener measure on $\cW$,
{\rm (iii)} $\cH=\cH^d$ is the $d$-dimensional Cameron-Martin space.
We denote by $(w_t)_{0\le t \le 1}$ the canonical realization 
of $d$-dimensional Brownian motion (i.e. the coordinate process).

Now we introduce a simple orthonormal basis of $\cH$.
First, set $\psi^{0,1}_t \equiv 1$.
For $n \ge 1$ and $1 \le m \le 2^{n-1}$,  set
\[
\psi^{n, m}_t = \left\{
\begin{array}{ll}
2^{(n-1)/2}, &  t \in [ (2m-2)2^{-n}, (2m-1)2^{-n}), \\
-2^{(n-1)/2}, &  t \in [ (2m-1)2^{-n}, 2m2^{-n}),   \\
0, &  \mbox{otherwise.}
\end{array}
\right.
\]
Denote by $\{ \mathbf{e}_i \}_{i=1}^d$ the canonical 
orthonormal basis of $\mathbb{R}^d$.
Then, it is well-known that 
\[
\{ \psi^{n, m}  \mathbf{e}_i \colon 
 n\ge 0,  \, 1 \le m \le 2^{n-1}\vee 1,   \, 1\le i \le d\}
\]
forms an orthonormal basis of $L^2 ([0,1], \mathbb{R}^d)$.
Since $L^2 ([0,1], \mathbb{R}^d)$ and $\cH$ are 
unitarily isometric, 
\[
\{ \varphi^{n, m}  \mathbf{e}_i \colon 
 n\ge 0,  \, 1 \le m \le 2^{n-1}\vee 1,   \, 1\le i \le d\}
\]
forms an orthonormal basis of $\cH$, where we set 
$ \varphi^{n, m}_t := \int_0^t  \psi^{n, m}_s ds$.

If we set $K_n$, $n\ge 1$, to be the linear span of 
\[
 \{ \varphi^{l, m}  \mathbf{e}_i \colon 
0\le l \le n-1,  \, 1 \le m \le 2^{l-1}\vee 1,   \, 1\le i \le d\},
\]
then 
$\cK=\{K_n\}^{\infty}_{n=1}$ is a non-decreasing, countable exhaustion of $\cH$ by finite dimensional subspaces.
Moreover, it is a routine to check that $P_{n} (h)= h(n)$ and 
$\bar{P}_{n} (w) = w(n)$ for all $n \ge 1$, $h\in \cH$ and $w\in \cW$.
Hence, we may apply Theorem \ref{aks>0} to 
the dyadic piecewise linear approximations of Brownian motion.
Finally, we remark that $\lim_{n\to\infty}\|w(n) -w\|_{\infty} =0$
for all $w \in \cW$
and $\lim_{n\to\infty}\|h(n) -h\|_{\cH} =0$ for all $h \in \cH$.
\end{example}

\section{Preliminaries from rough path theory}
\label{sec.RP}

In this section we recall the geometric rough path space 
with the H\"older or  Besov norm 
and quasi-sure properties of the rough path lift.
For basic properties of geometric rough path space with the H\"older 
topology, 
we refer to \cite{lcl, fvbk}.
For the geometric rough path space with the Besov topology, 
we refer to \cite[Appendix A.2]{fvbk}.
The quasi-sure properties of the 
rough path lift  are summarized in \cite{in1}.
From now on, $(\cW, \cH, \mu)$ stands for the $d$-dimensional 
classical Wiener space as in Example \ref{exm.DPL}.


In the first half of this section,
 we discuss deterministic aspects of rough path theory.
First, we work in the $\alpha$-H\"older rough path topology
with $\alpha \in (1/3, 1/2)$.
We consider an RDE with drift driven by 
${\bf w} \in G\Omega^{\textrm{H}}_{\alpha} ( {\mathbb R}^{d})$.
For vector fields $V_{i}: {\mathbb R}^e \to {\mathbb R}^e$ ($0 \le i \le d$), we consider the following RDE:
\begin{equation}\label{rde_x.def}
dx_t = \sum_{i=1}^d  V_i ( x_t) dw_t^i + V_0 ( x_t) dt,
\qquad
 x_0 =a \in \R^e.
\end{equation}
We assume that $V_i$, $0 \le i \le d$, is (at least) of $C_b^3$, 
that is, when viewed as an $\R^e$-valued function,
$V_i \in C^3_b (\R^e, \R^e)$.
It is then known that a unique global solution of \eqref{rde_x.def}
exists for every ${\bf w}$ and $a$. 
Moreover, Lyons' continuity theorem holds, that is, the map
\[
\Phi \colon G\Omega^{\textrm{H}}_{\alpha} ( {\mathbb R}^d)
\to 
\cC_a^{0, \alpha\textrm{-H}}({\mathbb R}^e)
\]
defined by $\Phi ({\bf w}) =x$ is locally Lipschitz continuous.
This map is called the Lyons-It\^o map.

\begin{remark}
We only study the first level paths of solutions of RDEs. 
Therefore, Lyons-It\^o map takes its values in a usual path space
and any formulation of RDEs will do.
\end{remark}

%

We introduce the skeleton ODE associated with 
RDE \eqref{rde_x.def}  and SDE \eqref{sde_x.def} below.
For $h \in \cH$, we consider the following ODE in the usual sense:
\begin{equation}\label{skele_x.def}
dx_t = \sum_{i=1}^d  V_i ( x_t) dh_t^i + V_0 ( x_t) dt,
\qquad
 x_0 =a \in \R^e.
\end{equation}
If $V_i$'s are of $C_b^{1}$, then a unique solution $x$ exists,
which is denoted by $\Psi (h)$.
Under the same condition, 
$\Psi\colon \cH\to \cC_a^{0, 1/2\textrm{-H}}({\mathbb R}^e)$ is locally Lipschitz continuous.
It should be noted that $\Psi (h) = \Phi ({\cal L}(h))$ (if $V_i$'s are of $C_b^3$).

%

Next we discuss Besov-type norms for rough paths.
We always assume that 
the Besov parameter $(\alpha, 4m)$ satisfy  the following 
conditions:
\begin{equation}\label{eq.amam}
\frac13 <\alpha < \frac12, \quad m \in \N, \quad
\al - \frac{1}{4m} > \frac13, \quad
 4m (\frac12 -\alpha)   >1.
\end{equation}
Observe that,
if the integer $m$ is chosen large enough for a given $\alpha \in (1/3, 1/2)$, then the two other inequalities in \eqref{eq.amam} are satisfied.
Heuristically, $\alpha$ plays a similar role to the H\"older parameter 
(see the Besov-H\"older embedding theorem below) and the auxiliary parameter
$4m$ is a very large even integer.

For $(\al, 4m)$ satisfying \eqref{eq.amam}, 
$G\Omega^{\textrm{B}}_{\alpha, 4m} ( {\mathbb R}^d)$ denotes the geometric rough path space 
over ${\mathbb R}^d$ with the $(\al, 4m)$-Besov norm.
It is defined to be the closure of 
$\{ \cL (k) \colon ~k\in \cC_0^{1\textrm{-H}} (\R^d) \}$ with respect to the
$(\al, 4m)$-Besov distance. The distance is given by 
\begin{align}
d_{\alpha, 4m}({\bf w}, \hat{\bf w}) 
&= \| {\bf w}^1- \hat{\bf w}^1 \|_{\al, 4m{\textrm{-B}}} 
+\| {\bf w}^2- \hat{\bf w}^2 \|_{2\al, 2m{\textrm{-B}}}
\nn\\
&
:=
\Bigl(
\iint_{0 \le s <t \le 1}  \frac{ | {\bf w}^1_{s,t}- \hat{\bf w}^1_{s,t}|^{4m}}
{|t-s|^{1 +4m\al }} 
dsdt
\Bigr)^{\tfrac{1}{4m}}
+
\Bigl(
\iint_{0 \le s <t \le 1}  \frac{ | {\bf w}^2_{s,t}- \hat{\bf w}^2_{s,t}|^{2m}}
{|t-s|^{1 +4m\al }} 
dsdt
\Bigr)^{\tfrac{1}{2m}}.
\nn
\end{align}
The homogeneous norm is denoted by 
$\vertiii{{\bf w}}_{\al, 4m{\textrm{-B}}} :=\| {\bf w}^1\|_{\al, 4m{\textrm{-B}}} +\| {\bf w}^2\|_{2\al, 2m{\textrm{-B}}}^{1/2}$.
It is known that $\{ \cL (h) \colon h\in \cH\}$ is dense in $G\Omega^{\textrm{B}}_{\alpha, 4m} ( {\mathbb R}^d)$.
By the Besov-H\"older embedding theorem for rough path spaces,
there is a continuous embedding $G\Omega^{\textrm{B}}_{\alpha, 4m} ( {\mathbb R}^d)
 \hookrightarrow G\Omega^{\textrm{H}}_{\alpha -(1/4m)} ( {\mathbb R}^d)$.
If $\al < \al' <1/2$, there is a continuous embedding
$G\Omega^{\textrm{H}}_{\alpha'} ( {\mathbb R}^d)
\hookrightarrow G\Omega^{\textrm{B}}_{\alpha, 4m} ( {\mathbb R}^d)$.
Basically, we will not write the first embedding explicitly.
(For example, if we write $\Phi({\bf w})$ for  
${\bf w} \in G\Omega^{\textrm{B}}_{\alpha, 4m} ( {\mathbb R}^d)$,
then it is actually the composition of the first embedding map above and $\Phi$ with respect to the
$\{\al -1/(4m)\}$-H\"older topology.)
It is known that the Young translation by $h \in {\cal H}$ works 
well on $G\Omega^{\textrm{B}}_{\alpha, 4m} ( {\mathbb R}^d)$ under (\ref{eq.amam}).
The map $({\bf w}, h) \mapsto T_h ({\bf w})$ is locally Lipschitz
continuous from $G\Omega^{\textrm{B}}_{\alpha, 4m} ( {\mathbb R}^d) \times {\cal H}$
to $G\Omega^{\textrm{B}}_{\alpha, 4m} ( {\mathbb R}^d)$, 
where 
$T_h ({\bf w})$ is the Young translation of ${\bf w}$ by $h$
(see \cite[Lemma 5.1]{in1}).
Recall that $T_h ({\bf w})$ is defined by 
\[
T_h ({\bf w})^1_{s,t} = {\bf w}^1_{s,t} +{\bf h}^1_{s,t}
\quad\mbox{and} \quad
T_h ({\bf w})^2_{s,t} = {\bf w}^2_{s,t} +{\bf h}^2_{s,t}
+ \int_s^t  {\bf w}^1_{s,u} \otimes dh_u
+\int_s^t  {\bf h}^1_{s,u} \otimes d_u  {\bf w}^1_{s,u}
\]
for $0\le s \le t \le 1$.
(The third and the fourth terms makes sense 
as a Riemann-Stieltjes and a Young integral, respectively.)


From here we discuss probabilistic aspects.
Suppose that $V_i$'s are of $C_b^{3}$
and let the notation as in Example \ref{exm.DPL}.
If ${\bf W}$ is Brownian rough path, 
that is, the natural (Stratonovich)
lift of $d$-dimensional Brownian motion $(w_t)_{0\le t \le 1}$, 
that is, 
\[
{\bf W}^1_{s,t} = w_t -w_s 
\quad\mbox{and} \quad
{\bf W}^2_{s,t} =  \int_s^t (w_u -w_s) \otimes \circ dw_u,
\qquad  0\le s \le t \le 1.
\]
Then, the process $(\Phi({\bf W})_t)_{0\le t \le 1}$
coincides $\mu$-a.s. with the solution $(X_t)_{0\le t \le 1}$ of the corresponding Stratonovich-type SDEs in the usual sense:
\begin{equation}\label{sde_x.def}
dX_t = \sum_{i=1}^d  V_i ( X_t)\circ dw_t^i + V_0 ( X_t) dt,
\qquad
 X_0 =a \in \R^e.
\end{equation}
Here, $\circ dw_t$ stands for the Stratonovich-type 
stochastic integral.
(The coefficients in \eqref{rde_x.def},
\eqref{skele_x.def} and \eqref{sde_x.def} 
are the same vector fields.)


Now we review quasi-sure properties of rough path lift map 
${\bf L}$ from ${\cal W}$ to $G\Omega^{\textrm{B}}_{\alpha, 4m} ( {\mathbb R}^d)$.
For $k\in \N$ and $w \in {\cal W}$,
we denote by $w(k)$ the $k$th dyadic piecewise linear approximation of $w$ 
associated with the partition $\{ j2^{-k} \colon 0 \le j \le 2^k\}$ of $[0,1]$.
We denote the natural lift of $w(k)$ by ${\cal L} (w(k))$.

For $(\alpha, 4m)$ satisfying \eqref{eq.amam}, we set 
\begin{equation}\label{def_Ztilde}
{\cal Z}_{\al, 4m} := \bigl\{ w \in {\cal W} \colon
\mbox{ $\{ {\cal L} (w(k)) \}_{k=1}^{\infty}$ is Cauchy in $G\Omega^{\textrm{B}}_{\alpha, 4m} ( {\mathbb R}^d)$} 
\bigr\}.   
\end{equation}
We define ${\bf L}: {\cal W} \to G\Omega^{\textrm{B}}_{\alpha, 4m} ( {\mathbb R}^d)$
by ${\bf L} (w) = \lim_{m\to \infty} {\cal L} (w(k))$ if $w \in {\cal Z}_{\al, 4m}$
and  we define ${\bf L} (w)={\bf 0}$ (the zero rough path)
 if $w \notin {\cal Z}_{\al, 4m}$.
 It is well-known that $\mu ({\cal Z}_{\al, 4m}) =0$.
Obviously, $w\mapsto  {\bf L} (w)$ is an 
everywhere-defined Borel measurable version of 
Brownian rough path ${\bf W}$ with respect to $\mu$.
(In what follows, when we write ${\bf W}$, 
it means this version.) 
 
It is easy to see that $\cH \subset {\cal Z}_{\al, 4m}$
and ${\bf L} (h) = {\cal L} (h)$ for all $h\in \cH$.
The scalar multiplication (i.e. the dilation)
and the Cameron-Martin translation leave
 ${\cal Z}_{\al, 4m}$ invariant.
Moreover, $c {\bf L} (w) = {\bf L} (cw)$ and $T_h({\bf L} (w))= {\bf L} (w+h)$ 
for all $w \in {\cal Z}_{\al, 4m}$, $c \in \R$, and $h \in {\cal H}$. 
It is known that ${\cal Z}_{\al, 4m}^c$ is slim, that is,
 the $(p,r)$-capacity of 
this set is zero for any $p \in (1,\infty)$ and $r \in \N$.
Therefore, from a viewpoint of quasi-sure analysis, 
the lift map ${\cal L}$ is well-defined.
Moreover, 
the map ${\cal W} \ni w \mapsto {\bf L} (w) \in G\Omega^{\textrm{B}}_{\alpha, 4m} ( {\mathbb R}^d) $
is $\infty$-quasi-continuous.
(This kind of $\infty$-quasi-continuity was first shown in \cite{aida}.)
Then, it immediately follows that the map
\[
\cW \ni \,\,
w \mapsto \Phi({\bf L} (w))  \,\, 
\in \cC_a^{0, \alpha-1/(4m)\textrm{-H}}({\mathbb R}^e)
\]
is an $\infty$-quasi-continuous version of $w \mapsto X$, 
where $X=(X_t)_{0\le t \le 1}$ is the solution of SDE \eqref{sde_x.def}
viewed as a path space-valued random variable.
%
%
%
%
\begin{remark}
The situation  described above can be summarized by the 
following commutative diagram:
\[
  \xymatrix{
  &  & G\Omega^{\textrm{B}}_{\alpha, 4m} ({\mathbb R}^d)   \ar[rd]^{\Phi} & \\
   & \cH \ar[ru]^{\cL}\ar@{^{(}-_>}[r]_{{\bf Incl}}
       & {\mathcal W}   \ar[r]_{{\bf Ito}\qquad } \ar[u]_{{\bf L}}   
             &\quad  \cC_a^{0, \alpha- (1/4m)\textrm{-H}}({\mathbb R}^e)
  }
\]
Here, ${\bf Incl}$ is the inclusion and ${\bf Ito}$ is
the usual It\^o map associated with SDE \eqref{sde_x.def}.
All maps above except ${\bf L}$ and ${\bf Ito}$ are continuous.
Note also that $\Psi = \Phi \circ \cL$.
\end{remark}

\begin{remark}
The first paper that used quasi-sure analysis 
for Brownian rough path is \cite{in0}.
In that paper, however, the rough path topology is the $p$-variation
topology with $2< p<3$.
The foundation of quasi-sure analysis 
for Brownian rough path in Besov or H\"older topology 
was laid by \cite{in1,aida}.
It was used for large deviations for pinned diffusion measures
in \cite{in1, in3, in4}.
A quasi-sure refinement of non-degeneracy property 
of Brownian signature was proved in \cite{bglq}.
It should also be noted that quasi-sure analysis for 
fractional Brownian rough path was studied in \cite{bgq,orv}.
\end{remark}


Before closing this section,  
let us recall Karhunen-Lo\'eve approximation, which will play an
important role in proofs of $\cK$-regularity 
and $\cK$-differentiability in the next section.
Fortunately, the dyadic piecewise linear approximation is also 
a Karhunen-Lo\'eve approximation
since $\bar{P}_{K_k} (w) = w(k)$ (see Example \ref{exm.DPL}).
It is easy to see that, for each fixed $k$, 
\begin{equation}\label{eq.0905-1}
T_{-w(k)} {\bf W} 
=\lim_{l\to\infty} T_{-w(k)}{\cal L} (w(l))
=\lim_{l\to\infty} {\cal L} (w(l) -w(k)),
\qquad
w\in {\cal Z}_{\al, 4m}.
\end{equation}
As one can naturally expects, the above quantity converges to 
the zero rough path as $k\to\infty$.
\begin{proposition}\label{pr.KL}
Let the notation be as above. Then, we have the following:
\begin{enumerate}
\item[{\rm (1)}]
There exists a positive constant $\eta$ independent of $k$
such that
\[
{\mathbb E} \bigl[ \exp \bigl( \eta 
\vertiii{ {\bf L} (w)}_{\al, 4m{\textrm{-B}}}^2
\bigr)\bigr]
\vee
\sup_{k \ge 1} {\mathbb E} \bigl[ \exp \bigl( \eta 
\vertiii{ \cL (w (k))}_{\al, 4m{\textrm{-B}}}^2
\bigr)
\bigr] <\infty.
\]
\item[{\rm (2)}]
For every $r \in [1,\infty)$ and $i =1, 2$,
$\lim_{k\to \infty}\| {\cal L} (w(k))^i - {\bf L}(w)^i\|_{i\al, 4m/i{\textrm{-B}}}  =0$
in $L^r(\mu)$.
\item[{\rm (3)}]
For every $r \in [1,\infty)$,
$\lim_{k\to \infty}\vertiii{T_{-w(k)} {\bf L} (w)}_{\al, 4m{\textrm{-B}}}=0$ in $L^r(\mu)$.
\end{enumerate}
\end{proposition}

\begin{proof}
If the rough path topology is $\beta$-H\"older with $1/3<\beta<1/2$,
these statements are proved in \cite[Theorem 15.47]{fvbk}.
Using the Besov-H\"older embedding theorem, we can easily prove this 
proposition, too.
\end{proof}


\section{$\cK$-differentiability of Lyons-It\^o map}
\label{sec.cK_diff}

In this section we show that 
the rough path lift map is uniformly $\cK$-regular
and Lyons-It\^o map is twice $\cK$-regularly differentiable. 
Technically, this section is the core of this paper.
These properties were already proved in \cite{ip}
for Gaussian rough paths with respect to the $p$-variation topology
under the condition called
the complementary Young regularity.
In this section, we will show these properties  for Brownian rough path
with respect to the Besov rough path topology and
also clean up arguments in \cite{ip}.
We keep the same notation as before.
Let $\cK=\{K_n\}_{n=1}^\infty$ be as in Example \ref{exm.DPL}.
We write
$\bar{P}_n (w) = w(n)$ and $P_n (h) = h(n)$ for $w \in \cW$
and $h\in \cH$.
We continue to assume \eqref{eq.amam} for the Besov parameter $(\alpha, 4m)$.
(For the rest of this paper, we study this particular exhaustion only.
We do not know what happens for a general exhaustion.)

First, we prove that
the rough path lift map is $\cK$-regular
and so is the solution of an RDE driven Brownian rough path.
Note that $\Phi \circ {\bf L}$ equals $\mu$-a.s. to the solution of 
RDE \eqref{rde_x.def} driven by ${\bf W}={\bf L}(w)$.

\begin{proposition}\label{pr.Kreg}
Let the notation be as above. Then, we have the following:
\begin{itemize}
\item[{\rm (1)}]
The measurable map ${\bf L}\colon \cW \to
G\Omega^{\textrm{B}}_{\alpha, 4m} ({\mathbb R}^d)$ is uniformly $\cK$-regular with $\cL$ as its regularization.
\item[{\rm (2)}]
Let $E_0$ be a Polish space and $\Lambda\colon G\Omega^{\textrm{B}}_{\alpha, 4m} ({\mathbb R}^d) \to E_0$ is locally Lipschitz continuous.
Then, $\Lambda \circ {\bf L}\colon \cW \to E_0$ is uniformly $\cK$-regular 
with $\Lambda \circ \cL$ as its regularization.
\item[{\rm (3)}]
If, in addition, 
$V_i$ is of $C_b^{3}$ for all $0 \le i \le d$, then 
$\Phi \circ 
{\bf L}\colon \cW \to \cC^{0, \alpha- (1/4m)\textrm{-H}}({\mathbb R}^e)$ is uniformly
$\cK$-regular with $\Psi$ as its regularization.
\end{itemize}
\end{proposition}


For the rest of this section we use the following notation.
We write $\cA :={\cal Z}_{\al, 4m}$,
which was defined by \eqref{def_Ztilde}.
It is important that this set is of full $\mu$-measure and
invariant under the translation by $h \in \cH$.
Write ${\bf W}^{*n} :=T_{-w(n)} {\bf L}(w)$.
By Proposition \ref{pr.KL}, 
 $\lim_{n\to\infty}{\bf W}^{*n} = {\bf 0}$ in probability.
We set $E^\prime =G\Omega^{\textrm{B}}_{\alpha, 4m} ({\mathbb R}^d)$, 
$G={\mathbf L}\colon \cW \to E^\prime$ and
$\tilde{G} = \cL\colon \cH \to E^\prime$. 
Similarly, we set 
$F=\Lambda \circ {\mathbf L}\colon \cW \to E_0$ and
$\tilde{F} =\Lambda \circ \cL\colon \cH \to E_0$. 
We will write 
$E = \cC^{0, \alpha- (1/4m)\textrm{-H}}({\mathbb R}^e)$.

\begin{proof}[Proof of Proposition \ref{pr.Kreg}]
In this proof, $\epsilon >0$ is arbitrary.
Set 
$G_n\colon\cW\times K_n \longmapsto E^\prime$ by 
$G_n (w,k) := T_k {\mathbf L}(w)$ if $w \in \cA$
and $G_n (w,k) = {\mathbf 0}$ if $w \notin \cA$.
Then, for all 
$w \in \cA$ and $k \in K_n$, we have 
\begin{align*}
G \circ \tau_k (w) &= {\mathbf L}(w+k) = T_k  {\mathbf L}(w) = G_n (w,k),
\\
G_n^\perp (w,k) &:= G_n(w,-\bar P_n (w)+k)
={\mathbf L}(w-\bar P_n (w)+k )
={\mathbf L}(\bar P_n^\perp (w)+k ) = T_k  {\mathbf W}^{*n}.
\end{align*}
Thanks to these explicit expressions, 
{\rm (i)} $\cK$-continuity of $G$ is now clear and 
{\rm (ii)}  we may and will view
 $G_n$ and $G_n^\perp$ as
  maps from  $\cW\times \cH$ to $E^\prime$.
(Then, they are actually independent of $n$.
Note that $G_n (w,k) = {\mathbf 0}=G_n^\perp (w,k)$ whenever $w \notin \cA$.)

We will check \eqref{rhoe2}. Take $w \in \cA$.
Note that $\tilde{G}(\bar{P}_n(w)+k) = \cL (w(n)+k )
= T_k \cL ( w (n))$ and 
$G_{n}(w,k) = T_k {\mathbf L}(w)$.
Since $\cL ( w (n)) \to  {\mathbf L}(w)$ as $n \to \infty$,
$\{ \cL ( w (n)) \}_{n=1}^\infty$ is bounded in $E^\prime$.
Since $T\colon E^\prime \times \cH \to E^\prime$ is locally Lipschitz continuous, we see that 
\begin{align*}
\sup_{\|k\|_{\cH}\leq r}
\rho_{E^\prime} (\tilde{G}(\bar{P}_n(w)+k),G_{n}(w,k))
&=
\sup_{\|k\|_{\cH}\leq r}
\rho_{E^\prime} (T_k \cL ( w (n)), T_k {\mathbf L}(w))
\\
&\le C_{r, w} \rho_{E^\prime} (\cL ( w (n)), {\mathbf L}(w))
\to 0
\quad
\mbox{as $n \to \infty$.}
\end{align*}
Here, 
$C_{r, w}$ is a positive constant which depends only on $r>0$
and $w \in \cA$ (and may vary from line to line).
Then, it immediately follows that, for every $m\in \N$, $\epsilon >0$ and $r>0$,
\[
 \lim_{n \to \infty} 
 \mu \Bigl(
\sup_{k\in K_m, \|k\|_{\cH}\le r}\rho_{E^\prime}
(\tilde{G}(\bar{P}_n(w)+k), G_{n}(w,k)) 
 \ge \epsilon\Bigr)=0.
 \]

Similarly, if $w\in \cA$ and $\vertiii{{\mathbf W}^{*n}} \le 1$, then we have
\begin{align}
\sup_{\|k\|_{\cH}\leq r}
\rho_{E^\prime} (\tilde{G}(h+k),G_n^{\perp}(w, P_n (h)+k))
&=
\sup_{\|k\|_{\cH}\leq r}
\rho_{E^\prime} ( T_{k+h} {\mathbf 0},  T_{k+ P_n (h)} {\mathbf W}^{*n})
\nn\\
&\le 
C_{r, h} \{ \rho_{E^\prime} ({\mathbf W}^{*n}, {\mathbf 0}) 
+ \|P_n (h)-h \|_{\cH} \}
\nn\\
&\le 
C_{r, h} \{ \vertiii{{\mathbf W}^{*n}} + \vertiii{{\mathbf W}^{*n}}^2
+ \| h(n) -h \|_{\cH} \}.
\label{eq.W^*}
\end{align}
Here, $C_{r, h}$ is a positive constant 
which depends only on $r>0$, $h\in \cH$. 
We can easily see from this that 
\begin{align} 
\lefteqn{
\mu \Bigl(
\sup_{k\in K_m, \|k\|_{\cH}\leq r}
\rho_{E^\prime}(\tilde{G}(h+k),  G^{\perp}_{n}(w, P_n(h)+k)) 
\ge \epsilon \Bigr)
}\nn\\
&\le 
\mu \bigl( \vertiii{{\mathbf W}^{*n}} \ge \tfrac{\epsilon}{3C_{r, h} }\bigr) 
+\mu \bigl( \vertiii{{\mathbf W}^{*n}}^2 \ge \tfrac{\epsilon}{3C_{r, h} }\bigr) 
 + \mu \bigl( \|h(n) -h \|_{\cH} \ge \tfrac{\epsilon}{3C_{r, h} }\bigr) 
 + \mu ( \vertiii{{\mathbf W}^{*n}} > 1 ).
 \nn
\end{align}
The right hand side tends to zero as $n\to\infty$ for every $m\in \N$, $\epsilon >0$, $h\in \cH$ and $r>0$.
Thus, we have shown (1).

Next we show (2).
Set also 
$F_n\colon\cW\times K_n\longmapsto E$ by 
$F_n (w,k) = \Lambda (G_n (w,k))$.
Take any $w \in \cA$ and $k \in \cH$. 
It is clear that $F \circ \tau_k (w)= F_n (w,k)$.
We also have $F_n^\perp (w,k) =  \Lambda (T_k  {\mathbf W}^{*n})$.
Again, we may and will view 
 $F_n$ and $F_n^\perp$ as maps from  $\cW\times \cH$ to $E_0$.
(Then, they are actually independent of $n$, too.)

Since $\Lambda$ is locally Lipschitz continuous
and both
$\tilde{G}(\bar{P}_n(w)+k)$ and $G_{n}(w,k)$ 
stay bounded as $n \in \N$ and $k\in \cH$ 
(with $\|k\|_{\cH} \le r$) vary,  we have 
\[
\sup_{ \|k\|_{\cH}\leq r}
\rho_{E_0} (\tilde{F}(\bar{P}_n(w)+k),F_{n}(w,k))
\le 
C_{r,w}
\sup_{ \|k\|_{\cH}\leq r}
\rho_{E^\prime} (\tilde{G}(\bar{P}_n(w)+k),G_{n}(w,k)).
\]
As we have seen, the right hand tends to zero as $n\to\infty$.
This implies that 
\[
 \lim_{n \to \infty} 
 \mu \Bigl(
\sup_{k\in K_m, \|k\|_{\cH}\le r}\rho_{E_0}
(\tilde{F}(\bar{P}_n(w)+k), F_{n}(w,k)) 
 \ge \epsilon\Bigr)=0
 \]
for every $m\in \N$, $\epsilon >0$ and $r>0$.


If $w\in \cA$ and $\vertiii{{\mathbf W}^{*n}} \le 1$,
we see from the local Lipschitz continuity of $\Lambda$ that 
\begin{align*}
\sup_{\|k\|_{\cH}\leq r}
\rho_{E_0} \bigl(\tilde{F}(h+k),F_n^{\perp}(w, P_n (h)+k) \bigr) 
&=
\sup_{\|k\|_{\cH}\leq r}
\rho_{E_0} \bigl(\Lambda (\tilde{G}(h+k)),  \, \Lambda (G_n^{\perp}(w, P_n (h)+k)) \bigr)
\nn\\
&=
\sup_{\|k\|_{\cH}\leq r}
\rho_{E_0} \bigl( \Lambda (T_{k+h} {\mathbf 0}),  \Lambda ( T_{k+ P_n (h)} {\mathbf W}^{*n}) \bigr)
\nn\\
&=
C_{r,h} \sup_{\|k\|_{\cH}\leq r}
\rho_{E^\prime} \bigl( T_{k+h} {\mathbf 0},   T_{k+ P_n (h)} {\mathbf W}^{*n} \bigr).
\end{align*}
where $C_{r, h}$ is a positive constant 
which depends only on $r>0$, $h\in \cH$. 
Recall that we have already computed the right hand.
So, we can show that
\[
\lim_{n\to \infty}  
\mu \Bigl(
\sup_{k\in K_m, \|k\|_{\cH}\leq r}
\rho_{E_0}(\tilde{F}(h+k),  F^{\perp}_{n}(w, P_n(h)+k)) 
\ge \epsilon \Bigr) =0
\]
for every $m\in \N$, $\epsilon >0$, $h\in \cH$ and $r>0$
in exactly the same way as above. Thus, we have shown (2).

Finally, (3) is just a special case of (2) since $\Phi$ is locally Lipschitz
continuous.
Note that $\tilde{F}=\Phi \circ \cL=\Psi\colon \cH \to E$ in this case, which is the solution map of the 
skeleton ODE \eqref{skele_x.def}.
\end{proof}



We consider derivatives of the 
solution map $\Psi\colon \cH\to \cC_a^{0, 1/2\textrm{-H}}({\mathbb R}^e)
\subset \cC^{0,1/2\textrm{-H}}({\mathbb R}^e)$
of the skeleton ODE \eqref{skele_x.def}.
For brevity, we write $\sigma = [V_1, \ldots, V_d]$ and $b =V_0$ 
and view them as an $e \times d$ matrix-valued 
and an ${\mathbb R}^e$-valued function, respectively.
In what follows, we assume these coefficients are of $C^5_b$
for simplicity.
Then, \eqref{skele_x.def} simply reads
\[
dx_t = \sigma (x_t) dh_t + b (x_t) dt, \qquad  x_0 =a \in \R^e.
\]
It is well-known that $\Psi$ (i.e. $h \mapsto x=x(h)$)
is Fr\'echet-$C^2$.
Moreover, its directional derivatives satisfy 
a simple ODE, which can be obtained as a formal 
differentiation of the above ODE.
For example, consider 
$D_l x_t$ and $D^2_{l, \tilde{l}} x_t$ for $l \in \cH$, 
where $D$ stands for the Fr\'echet derivative on $\cH$ and 
$l\in\cH$ is a direction of differentiation.
If those are denoted by $\xi_t^{[1]}= \xi_t^{[1]} (h;l)$ and
 $\xi_t^{[2]}=\xi_t^{[2]} (h; l, \tilde{l})$,  their ODEs explicitly read 
\begin{equation}
d \xi_t^{[1]}
= 
\nabla \sigma (x_t ) \la  \xi_t^{[1]}, dh_t \ra
+ 
\nabla b(x_t ) \la  \xi_t^{[1]}\ra dt
+
\sigma (x_t ) dl_t,
\qquad
\xi_0^{[1]}  =0\in {\mathbb R}^e,
\label{heu_1.eq}
\end{equation}
and 
\begin{align}
d \xi_t^{[2]}
&= 
\nabla \sigma (x_t) \la \xi_t^{[2]}, dh_t \ra 
+
\nabla b (x_t) \la \xi_t^{[2]} \ra  dt
\nn\\
&\qquad
+
\nabla^2 \sigma (x_t) \la \xi_t^{[1]}(h;l),  \xi_t^{[1]}(h; \tilde{l}), dh_t\ra
+
\nabla \sigma (x_t) \la \xi_t^{[1]}(h; \tilde{l}),dl_t\ra 
\nn\\
&\qquad
+
\nabla \sigma (x_t) \la \xi_t^{[1]}(h;l),d\tilde{l}_t\ra 
+
\nabla^2 b (x_t) \la \xi_t^{[1]}(h;l),  \xi_t^{[1]}(h; \tilde{l})\ra dt,
\qquad
\xi_0^{[2]}=0\in {\mathbb R}^e,
\label{heu_2.eq}
\end{align}
respectively.
Note that both are simple first-order ODEs and therefore 
can be solved by the variation of constants formula
for every $h, l, \tilde{l}\in \cH$.
It is also standard to show that the map
\[
\cH \times \cH \times \cH \ni
(h,l, \tilde{l}) \mapsto (\Psi(h), \xi^{[1]}(h;l), \xi^{[2]}(h;l, \tilde{l}))
\in
\cC^{0, 1/2\textrm{-H}} (({\mathbb R}^e)^{\oplus 3} 
)
\]
is locally Lipschitz continuous.

Now we get back to RDEs. 
The RDE driven by ${\bf w}\in G\Omega^{\textrm{B}}_{\alpha, 4m} ({\mathbb R}^d)$ 
and $l\in \cH$
which correspond to \eqref{heu_1.eq}
is given as follows:
\begin{align}
d \xi_t^{[1]}
&= 
\nabla \sigma (x_t ) \la  \xi_t^{[1]}, dw_t \ra
+ 
\nabla b(x_t ) \la  \xi_t^{[1]}\ra dt
+
\sigma (x_t ) dl_t,
\qquad
\xi_0^{[1]}  =0\in {\mathbb R}^e.
\label{heu_3.eq}
\end{align}
We write $\xi_t^{[1]}=  \xi_t^{[1]} ({\bf w}; l)$ when necessary.
Likewise, the RDE driven by ${\bf w}\in G\Omega^{\textrm{B}}_{\alpha, 4m} ({\mathbb R}^d)$  and $l, \tilde{l} \in \cH$
which correspond to \eqref{heu_2.eq}
is given as follows:
\begin{align}
d \xi_t^{[2]}
&= 
\nabla \sigma (x_t) \la \xi_t^{[2]}, dw_t \ra 
+
\nabla b (x_t) \la \xi_t^{[2]} \ra  dt
\nn\\
&\quad
+
\nabla^2 \sigma (x_t) \la \xi_t^{[1]}({\bf w}; l),  \xi_t^{[1]}({\bf w}; \tilde{l}), dw_t\ra
+
\nabla \sigma (x_t) \la \xi_t^{[1]}({\bf w}; \tilde{l}),dl_t\ra 
\nn\\
&\quad
+
\nabla \sigma (x_t) \la \xi_t^{[1]}({\bf w}; l),d\tilde{l}_t\ra 
+
\nabla^2 b (x_t) \la \xi_t^{[1]}({\bf w}; l),  \xi_t^{[1]}({\bf w}; \tilde{l}) \ra dt,
\quad
\xi_0^{[2]}=0\in {\mathbb R}^e.
\label{heu_4.eq}
\end{align}
We write $\xi_t^{[2]}=  \xi_t^{[2]} ({\bf w}; l, \tilde{l})$ when necessary.
Since \eqref{heu_3.eq}--\eqref{heu_4.eq} are first-order RDEs, 
it is known that the system of three RDEs \eqref{rde_x.def},
\eqref{heu_3.eq} and \eqref{heu_4.eq}
has a unique global solution for every $({\bf w}, l)$.
Moreover, a rough path version of the variation of 
constants formula holds for $\xi^{[1]}$ and $\xi^{[2]}$, too.
(See \cite[Section 10.7]{fvbk}  and \cite{in2}  for example.)
If $h\in \cH$, we have
$\xi^{[1]} (h; l)= \xi^{[1]} ({\mathcal L} (h); l)$ and 
$\xi^{[2]} (h; l, \tilde{l})= \xi^{[2]} ({\mathcal L} (h); l, \tilde{l})$.
Since no explosion can happen, Lyons' continuity theorem 
still holds for this system of three RDEs. 
In particular,
the following map is locally Lipschitz continuous:
\begin{align}\label{eq.der.Lip.conti}
G\Omega^{\textrm{B}}_{\alpha, 4m} ({\mathbb R}^d)
\times \cH\times \cH   &\ni ({\bf w}, l, \tilde{l})
\nn\\
&\qquad \mapsto 
(\Phi({\bf w}), \xi^{[1]}({\bf w};l), \xi^{[2]}({\bf w}; l, \tilde{l}))
\in
\cC^{0, \alpha-(1/4m) \textrm{-H}}(({\mathbb R}^e)^{\oplus 3} ).
\end{align}
Here, $\Phi({\bf w})$ is the solution of RDE \eqref{rde_x.def}.
This property will play a key role.

\begin{lemma}\label{lem.fre_diff}
Let $F_n \colon \cW\times \cH \to E:=\cC^{0, \alpha- (1/4m)\textrm{-H}}({\mathbb R}^e)$ be as in the proof of Proposition \ref{pr.Kreg},
namely,
\[
F_n (w,k) =
\begin{cases}
\Phi(T_k {\bf L}(w)) & (\mbox{if $w \in \cA$ and $k\in \cH$}),
\\
\Phi({\bf 0}) & (\mbox{if $w \notin \cA$ and $k\in \cH$}).
\end{cases}
\]
Then, for each $w \in \cW$, $F_n (w, \bullet)\colon \cH \to E$ is of
Fr\'echet-$C^2$.
Moreover, for each $w \in \cA$, we have 
\[
D_l F_n (w,k)=\xi^{[1]}( T_k {\mathbf L}(w); l), 
\quad
D^2_{l, \tilde{l}} F_n (w,k)=\xi^{[2]}( T_k {\mathbf L}(w); l, \tilde{l}), 
\qquad
k, l, \tilde{l} \in \cH.
\]
Here, $D$ is the Fr\'echet derivative acting on the $k$-variable.
\end{lemma}

\begin{proof}
As is well-known, Fr\'echet-$C^j$ and G\^ateaux-$C^j$
are equivalent for all $j \ge 1$. 
So, we only consider G\^ateaux derivatives.
The case $w \notin \cA$ is obvious. So, we pick any $w \in \cA$
and will fix it in what follows.

First, we calculate the first-order derivative in the direction $l$.
For $m \in \N$, $w(m) \in \cH$ and therefore we already know that
$\lim_{m\to\infty} F_n (w(m),k)=F_n (w,k)$ and 
\[
D_l F_n (w(m),k)=\xi^{[1]}( T_k {\mathcal L}(w(m)); l).
\]
The right hand converges to 
$\xi^{[1]}( T_k {\mathbf L}(w); l)$ in $E$ as $m\to\infty$ uniformly on
every ball of $\cH$.
Indeed, if $\|k\|_{\cH} \le r$ and $\|l\|_{\cH} \le r'$ ($r, r' >0$), then 
by the local Lipschitz continuity in \eqref{eq.der.Lip.conti}
and that of $T$, we have
\[
\| \xi^{[1]}( T_k \cL ( w (m)); l)
- \xi^{[1]}(T_k{\mathbf L}(w); l) \|_E 
\le 
C_{r, r',w} \, \rho_{E^\prime} (\cL ( w (m)), {\mathbf L}(w))
\to 0 \quad \mbox{as $m\to \infty$},
\]
where $E^\prime :=G\Omega^{\textrm{B}}_{\alpha, 4m} ({\mathbb R}^d)$
and $C_{r, r',w}>0$ is constant which depends only on $r, r',w$.

This uniform convergence in $k$ (for a fixed $l$) yields
the desired formula for the first derivative.
This can be verified as follows.
For fixed $k$ and $l$,
 set $\chi_m (\tau) = F_n (w(m),k +\tau l)$, $\tau \in (-1,1)$.
Obviously, $\lim_{m\to\infty}\chi_m (\tau) 
=F_n (w,k+\tau l)$ for every $\tau$ and
$(\partial/ \partial \tau )\chi_m (\tau) = 
\xi^{[1]}( T_{k +\tau l} \cL ( w (m)); l)$. As we have seen, 
$(\partial/ \partial \tau )\chi_m (\tau)$ converges to 
$\xi^{[1]}(T_{k +\tau l} {\mathbf L}(w); l)$ uniformly in $\tau\in (-1,1)$.
Hence, we have $(\partial/ \partial \tau ) F_n (w,k+\tau l) = 
 \xi^{[1]}(T_{k +\tau l} {\mathbf L}(w); l)$.
  Setting $\tau =0$, we obtain the formula.

By a similar argument we can show the continuity of 
$k \mapsto DF_n(w,k)$ as follows: 
If $\|k\|_{\cH},  \|\tilde{k}\|_{\cH}\le r$  ($r>0$), then 
by the local Lipschitz continuity we have
\begin{align} 
\| DF_n(w,k) - DF_n(w,\tilde{k}) \|_{\cH \to E}
&=
\sup_{\|l\| \le 1} \| D_l F_n(w,k) - D_l F_n(w,\tilde{k}) \|_{E}
\nn\\
&=
\| \xi^{[1]}( T_k {\mathbf L}(w); l)
- \xi^{[1]}(T_{\tilde{k}} {\mathbf L}(w); l) \|_E 
\nn\\
&\le 
C_{r,w} \, \|k -\tilde{k}\|_{\cH},
\end{align}
where $\| \,\cdot\, \|_{\cH \to E}$ is  the 
operator norm for bounded operators from $\cH$ to $E$
and $C_{r,w}>0$ is constant which depends only on $r,w$.
Thus, we have seen that $F_n (w, \bullet)\colon \cH \to E$ is of $C^1$.

Starting with the known fact that
\[
D^2_{l,\tilde{l}} F_n (w(m),k)=D_{\tilde{l} }\xi^{[1]}( T_k {\mathcal L}(w(m)); l) 
= \xi^{[2]}( T_k {\mathcal L}(w(m)); l,\tilde{l}),
\]
we can calculate the second-order derivative, too.
Since the proof is essentially the same as in the first-order case, 
we omit it.
\end{proof}


The following is the main result of this section.
It immediately implies that It\^o map   
$w \mapsto \Phi ( {\bf L} (w))_t$ at a fixed time $t$
is also twice $\cK$-regularly differentiable.
\begin{proposition}\label{pr.Kdiff}
Let the notation be as above
and assume that $V_i$ is of $C_b^{5}$ for all $0 \le i \le d$.
Then, 
$\Phi \circ 
{\bf L}\colon \cW \to  E :=
\cC^{0, \alpha-(1/4m) \textrm{-H}}({\mathbb R}^e)$
is twice $\cK$-regularly differentiable with $\Psi$ as its regularization.
\end{proposition}

\begin{proof}
We use the same notation as in the proof of 
Proposition \ref{pr.Kreg} (2).
An unimportant positive constant which depends only 
on the parameter $\star$ is denoted by $C_\star$,
which may vary from line to line.

We will prove \eqref{rhoe3} for $l=2$ by estimating 
\begin{align}
\|\tilde{F}(\bar{P}_n(w)+\bullet)-F_{n}(w,\bullet))\|_{C_{b}^2(B_{\cH}(0,r), E)} 
\vee \|\tilde{F}(h+\bullet)-F^{\perp}_{n}(w, P_n(h)+\bullet)\|_{C_{b}^2(B_{\cH}(0,r), E)} 
\nn
 \end{align} 
 for every $w\in \cA$, $r>0$ and $h \in \cH$.
Convergence of the zeroth order was already shown 
in the proof of Proposition \ref{pr.Kreg} (2).

Now we calculate the first order derivatives.
For the rest of the proof,  let $r, r'>0$, $w \in {\mathcal A}$, 
$k, l, h \in {\mathcal H}$.
$D_l   F_{n}(w,\bullet)$ was calculated in Lemma \ref{lem.fre_diff} above.
Since $\tilde{F} = \Psi$, we have
\begin{align*}
D_l \tilde{F}(\bar{P}_n(w)+\bullet) \vert_{\bullet =k}
&=
\xi^{[1]}(w(n)+ k; l)
=
\xi^{[1]}( T_k \cL ( w (n)); l).
\end{align*}
Due to the local Lipschitz continuity of $\xi^{[1]}$ we mentioned in 
\eqref{eq.der.Lip.conti}, we have that, 
if $\|k\|_{\cH} \le r$, then we have
\begin{align}
\lefteqn{
\sup_{\|k\|_{\cH} \le r}
\|D \tilde{F}(\bar{P}_n(w)+\bullet)\vert_{\bullet =k} 
- D F_{n}(w,\bullet))\vert_{\bullet =k}
\|_{\cH \to E}
}
\nn\\
&=
\sup_{\|k\|_{\cH} \le r} \sup_{\|l\|_{\cH} \le 1}
 \bigl\|
D_l \tilde{F}(\bar{P}_n(w)+\bullet) \vert_{\bullet =k} 
- D_l F_{n}(w,\bullet))\vert_{\bullet =k}
 \bigr\|_E 
\nn\\
&=
\sup_{\|k\|_{\cH} \le r}  \sup_{\|l\|_{\cH} \le 1} 
   \| \xi^{[1]}( T_k \cL ( w (n)); l)
           - \xi^{[1]}(T_k{\mathbf L}(w); l) \|_E
\nn\\
&\le 
C_{r, w} \rho_{E^\prime} (\cL ( w (n)), {\mathbf L}(w))
\to 0
\quad
\mbox{as $n \to \infty$}
\nn
\end{align}
for every $w\in\cA$ and $r>0$. Then, it immediately follows that
\[
\lim_{n \to \infty} 
 \mu \Bigl(
\sup_{k\in K_m, \|k\|_{\cH}\le r}
\|D \tilde{F}(\bar{P}_n(w)+\bullet)\vert_{\bullet =k} 
- D F_{n}(w,\bullet))\vert_{\bullet =k}
\|_{\cH \to E} \ge \epsilon\Bigr)=0
\]
for every $m\in \N$, $\epsilon >0$ and $r>0$.

Since $F_n^\perp (w,k) =F_n(w,-w(n) +k)$,
we see that 
\[
D_l F_n^\perp (w,\bullet)\vert_{\bullet =k}
= \xi^{[1]}( T_{k -w(n)}  {\mathbf L}(w); l) =\xi^{[1]}( T_{k}  {\bf W}^{*n}; l).
\]
Hence, if $h\in \cH$, $w\in \cA$ and $\vertiii{{\mathbf W}^{*n}} \le 1$, then
\begin{align}\label{eq.200323-4}
\lefteqn{
\sup_{\|k\|_{\cH} \le r} 
\bigl\|
 D\tilde{F}(h+\bullet)\vert_{\bullet =k} 
- D F^{\perp}_{n}(w, P_n(h)+\bullet) \vert_{\bullet =k}   
\bigr\|_{\cH\to E }
}
\nn\\
&=
\sup_{\|k\|_{\cH} \le r} \sup_{\|l\|_{\cH} \le 1} 
\bigl\|
D_l  \tilde{F}(h+\bullet)\vert_{\bullet =k} 
- D_l F^{\perp}_{n}(w, P_n(h)+\bullet) \vert_{\bullet =k}   \bigr\|_E 
\nn\\
&\le
\sup_{\|k\|_{\cH} \le r} \sup_{\|l\|_{\cH} \le 1} 
\|  
\xi^{[1]} (T_{h+k} {\bf 0}; l) - \xi^{[1]}( T_{k+P_n(h)}  {\bf W}^{*n};  l)
\|_E
\nn\\
&\le 
C_{r, h} \{
\rho_{E^\prime} ( {\mathbf W}^{*n}, {\mathbf 0}) +\| h-h(n)\|_{\cH} \}
\nn
\end{align}
for every $r>0$ and $h \in\cH$.
Recall that the right hand side is essentially the same
as the right hand side of \eqref{eq.W^*} in the proof of Proposition \ref{pr.Kreg}. So, we can show in the same way that 
\[
\lim_{n\to\infty} \mu \Bigl(
\sup_{k\in K_m, \|k\|_{\cH}\le r}
\bigl\|
 D\tilde{F}(h+\bullet)\vert_{\bullet =k} 
- D F^{\perp}_{n}(w, P_n(h)+\bullet) \vert_{\bullet =k}   
\bigr\|_{\cH\to E }
\ge \epsilon \Bigr) =0
\]
for every $m\in \N$, $\epsilon >0$, $r>0$ and $h \in\cH$.

Next we calculate the second order derivatives.
We can easily see that
\begin{align*}
D^2_{l, \tilde{l}} \tilde{F}(\bar{P}_n(w)+\bullet) \vert_{\bullet =k}
&=
\xi^{[2]}(w(n)+ k; l, \tilde{l})
=
\xi^{[2]}( T_k \cL ( w (n)); l, \tilde{l}).
\end{align*}
Due to Lemma \ref{lem.fre_diff} and 
the local Lipschitz continuity of $\xi^{[2]}$ we mentioned in 
\eqref{eq.der.Lip.conti}, we have the following: 
\begin{align}\label{eq.200324-1}
\lefteqn{
\sup_{\|k\|_{\cH} \le r}
\|D^2 \tilde{F}(\bar{P}_n(w)+\bullet)\vert_{\bullet =k} 
- D^2 F_{n}(w,\bullet))\vert_{\bullet =k}
\|_{\cH \times \cH \to E}
}
\nn\\
&=
\sup_{\|k\|_{\cH} \le r}\sup_{\|l\|_{\cH} \vee \|\hat{l}\|_{\cH} \le 1}
 \| \xi^{[2]}( T_k \cL ( w (n)); l, \hat{l})
   - \xi^{[2]}( T_k {\mathbf L}(w); l, \hat{l}) \|_E
\nn\\
&\le 
C_{r, w} \rho_{E^\prime} (\cL ( w (n)), {\mathbf L}(w))
\to 0
\quad
\mbox{as $n \to \infty$}
\end{align}
for every $w \in \cA$ and $r>0$.
Here, $\| \,\cdot\, \|_{\cH \times \cH \to E}$ is the standard
norm for bounded bilinear maps from $\cH\times \cH $ to $E$.
From this we see that
\[
\lim_{n \to \infty} 
 \mu \Bigl(
\sup_{k\in K_m, \|k\|_{\cH}\le r}
\|D^2 \tilde{F}(\bar{P}_n(w)+\bullet)\vert_{\bullet =k} 
- D^2 F_{n}(w,\bullet))\vert_{\bullet =k}
\|_{\cH\times \cH \to E} \ge \epsilon\Bigr)=0
\]
for every $m\in \N$, $\epsilon >0$ and $r>0$.

In a similar way as above,  we see that 
$D^2_{l, \tilde{l}} F_n^\perp (w,\bullet)\vert_{\bullet =k}
=
\xi^{[2]}( T_{k}  {\bf W}^{*n}; l, \tilde{l})$.
Hence,  if $h\in \cH$, $w\in \cA$ and $\vertiii{{\mathbf W}^{*n}} \le 1$, 
then we have
\begin{align}\label{eq.200324-6}
\lefteqn{
\sup_{\|k\|_{\cH} \le r} \bigl\|
D^2 \tilde{F}(h+\bullet)\vert_{\bullet =k} 
- D^2 F^{\perp}_{n}(w, P_n(h)+\bullet) \vert_{\bullet =k}   
\bigr\|_{\cH\times\cH \to E} 
}
\nn\\
&=
\sup_{\|k\|_{\cH} \le r}\sup_{\|l\|_{\cH} \vee \|\hat{l}\|_{\cH} \le 1}
\bigl\|
D^2_{l,\hat{l}}
 \tilde{F}(h+\bullet)\vert_{\bullet =k} 
- D^2_{l,\hat{l}}
 F^{\perp}_{n}(w, P_n(h)+\bullet) \vert_{\bullet =k}   \bigr\|_E 
 \nn\\
&\le
\sup_{\|k\|_{\cH} \le r}\sup_{\|l\|_{\cH} \vee \|\hat{l}\|_{\cH} \le 1}
\|   \xi^{[2]} (T_{h+k} {\bf 0}; l, \hat{l}) - 
\xi^{[2]}( T_{k+P_n(h)}  {\bf W}^{*n}; l, \hat{l})
\|_E
\nn\\
&\le 
C_{r, h} \{
\rho_{E^\prime} ( {\mathbf W}^{*n}, {\mathbf 0})
+  \| h-h(n)\|_{\cH} \}
\nn
\end{align}
for every $r>0$.
As we have seen, this implies again that 
\[
\lim_{n\to\infty} \mu \Bigl(
\sup_{k\in K_m, \|k\|_{\cH}\le r}
\bigl\|
 D^2\tilde{F}(h+\bullet)\vert_{\bullet =k} 
- D^2 F^{\perp}_{n}(w, P_n(h)+\bullet) \vert_{\bullet =k}   
\bigr\|_{\cH\times \cH\to E }
\ge \epsilon \Bigr) =0
\]
for every $m\in \N$, $\epsilon >0$, $r>0$ and $h \in\cH$.
This completes the proof.
\end{proof}
%
%
%
%
%

\section{Support theorem on geometric rough path space}
\label{sec.SUP}

In this section we consider SDE \eqref{sde_x.def}
with $C^\infty_b$-vector fields $V_i~(0\le i \le d)$.
When we emphasize the starting point $a$, we write $X_t =X(t, a)$.
Similarly, the corresponding Lyons-It\^o map is denoted by 
$\Phi^a =\Phi$.
(Similarly, the deterministic It\^o map for the skeleton ODE 
is denoted by $\Psi^a =\Psi$.)
Recall that $\Phi^a ({\bf L}(w))$ is an $\infty$-quasi-continuous 
modification of $X(\,\cdot\,, a)$ for every $a$.
We continue to assume \eqref{eq.amam} for the Besov parameter $(\alpha, 4m)$.
The diffusion semigroup associated with this SDE is denoted 
by $(T_t)_{0\le t \le 1}$, that is, $T_t f (a) :=\E [X(t, a)]$
for every bounded continuous function $f\colon \R^e \to \R$.

Let $\cV$ be a linear subspace of $\R^e$ with dimension 
$e' ~(1 \le e' \le e)$.
The inner product of $\cV$ is a restriction of that of $\R^e$
and therefore the Lebesgue measure on $\R^e$ is uniquely determined.
The orthogonal projection from $\R^e$ onto $\cV$ is denoted by $\Pi$.
We are interested in the law of  $Y(\cdot, a) := \Pi (X(\cdot,a))$, in particular,  
when it is pinned at $b \in \cV$ at $t=1$.
It is well-known that $Y(t, a)$ is a
${\bf D}_\infty$-Wiener functional for every $t$ and $a$.

Suppose that Malliavin covariance of $Y(1,a)$ is non-degenerate.
Then, $\delta_b \circ Y(1,a) =\delta_b (Y(1,a)) \in \tilde{\bf D}_{-\infty}$
is well-defined as a positive Watanabe distribution.
By Sugita's theorem (see \cite{sugita} or \cite[p.101]{ma}), 
$\delta_b (Y(1,a))$ corresponds to
 a unique finite Borel measure on $\cW$, which is denoted by 
 $\hat\mu_{a,b}$.
The correspondence is given by 
\[
\E [G \,  \delta_b (Y(1,a))] = \int_{\cW}  \tilde{G}(w) \hat\mu_{a,b}(dw),
\qquad
G \in {\bf D}_{\infty},
\]
where $\tilde{G}$ is an $\infty$-quasi-continuous modification of $G$.
If the total mass $\E [ \delta_b (Y(1,a))]<\infty$ is positive, then 
$\mu_{a,b}:=\E [ \delta_b (Y(1,a))]^{-1} \hat\mu_{a,b}$ is a probability measure.
By Theorem \ref{aks>0} and Proposition \ref{pr.Kdiff},
$\E [ \delta_b (Y(1,a))] >0$ if and only if 
\begin{equation}\label{eq.non_emp}
\{ h \in \cH \colon   \mbox{$D\Pi \Psi^a(h)_1
\colon \cH\rightarrow \mathcal{V}$ is surjective, \,
$\Pi \Psi^a(h)_1=b$}
\} \neq \emptyset.
\end{equation}
(Here, we used Theorem \ref{aks>0} with $G\equiv 1$, $F(w) = Y(1,a) =\Pi \Phi^a(w)_1$ and $\tilde{F}(h) =\Pi \Psi^a(h)_1$.)
This measure  does not charge a slim set. 
So, its rough path lift ${\bf L}_* \mu_{a,b} = \mu_{a,b} \circ {\bf L}^{-1}$ is a Borel probability 
measure on $G\Omega^{\textrm{B}}_{\alpha, 4m} ( {\mathbb R}^d)$,
which will be denoted by $\nu_{a,b}$.

\begin{theorem} \label{or.supp.GRP}
Let the notation and the situation be as above.
Assume \eqref{eq.non_emp} and non-degeneracy of $Y(1,a)$.
Then, the support of $\nu_{a,b}$  equals the closure of 
\begin{equation}\label{eq.KLM}
\{ \cL (h)  \colon   \mbox{$h \in \cH$, $D\Pi \Psi^a(h)_1
\colon \cH\rightarrow \mathcal{V}$ is surjective, 
$\Pi \Psi^a(h)_1=b$}
\}
\end{equation}
in $G\Omega^{\textrm{B}}_{\alpha, 4m} ( {\mathbb R}^d)$.
\end{theorem}

\begin{proof} 
We use Aida-Kusuoka-Stroock's positivity theorem (Theorem \ref{aks>0}).
Let ${\bf z}=({\bf z}^1, {\bf z}^2)\in G\Omega^{\textrm{B}}_{\alpha, 4m} ( {\mathbb R}^d)$.
For $r>0$, we set
\[
B({\bf z}, r) :=
\bigl\{ {\bf w}\in G\Omega^{\textrm{B}}_{\alpha, 4m} ( {\mathbb R}^d)
\colon 
\| {\bf w}^1- {\bf z}^1 \|_{\al, 4m{\textrm{-B}}}^{4m}
+\| {\bf w}^2- {\bf z}^2 \|_{2\al, 2m{\textrm{-B}}}^{2m} <r^{4m}
\bigr\}.
\]
Then, $\{ B({\bf z}, r) \}_{r>0}$ forms a fundamental system of 
open neighborhood around ${\bf z}$.
Let $\chi\colon [0,\infty)\to [0,1]$ be a non-increasing $C^\infty$-function 
such that $\chi \equiv 1$ on $[0,1]$ and $\chi \equiv 0$ on $[2^{4m},\infty)$.
Set a non-negative ${\bf D}_\infty$-functional by
\[
G_{{\bf z}, r} (w) := \chi \Bigl( \frac{\| {\bf L}(w)^1- {\bf z}^1 \|_{\al, 4m{\textrm{-B}}}^{4m}
+\| {\bf L}(w)^2- {\bf z}^2 \|_{2\al, 2m{\textrm{-B}}}^{2m}}{r^{4m}}
\Bigr).
\]
It is easy to see from Proposition \ref{pr.Kreg}
that $G_{{\bf z}, r}$ is uniformly $\cK$-regular.
Obviously, 
\begin{equation}\label{eq.indG}
 {\bf 1}_{B({\bf z}, r)} ({\bf L}(w)) \le 
G_{{\bf z}, r} (w) \le {\bf 1}_{B({\bf z}, 2r)} ({\bf L}(w)), 
\qquad \mbox{for quasi-every $w \in \cW$.}
\end{equation}
As usual, ${\bf 1}_C$ denotes the indicator function of a subset $C$.

First, we show that $\bar{A} \subset {\rm supp} (\nu_{a,b})$,
where $A$ stands for the set in \eqref{eq.KLM}.
To do so, it suffices to show that $\nu_{a,b} (B(\cL (h), 2r)) >0$ for
every $\cL (h) \in A$  and $r>0$.
By \eqref{eq.indG} we have
\begin{align}
\nu_{a,b} (B(\cL (h), 2r)) 
&= \int_{\cW}  {\bf 1}_{B(\cL (h), 2r)} ({\bf L}(w))   \mu_{a,b}(dw)
\nn\\
&\ge 
\int_{\cW} G_{\cL (h), r} (w)   \mu_{a,b}(dw)
= \E [G_{\cL (h), r} \, \delta_b (Y(1,a))]/Z. 
\label{eq.posiposi}
\end{align}
Here, we used $\infty$-quasi-continuity of $G$ and wrote $Z:= \E [ \delta_b (Y(1,a))]>0$ for simplicity.
Recall that $Y(1,a)$ is non-degenerate by assumption and 
twice $\cK$-regularly differentiable 
(with its regularization $\Pi \Psi^a (\cdot)_1$) by Proposition \ref{pr.Kdiff}.
Therefore, we can use Theorem \ref{aks>0}  to the right hand side of \eqref{eq.posiposi} is positive for every $r>0$.
(Note that $h$ itself satisfies the second one of the two 
equivalent conditions in Theorem \ref{aks>0}.)

Next, we show that $\bar{A} \supset {\rm supp} (\nu_{a,b})$.
It is sufficient to show that for every ${\bf z}\notin \bar{A}$,
there exists $r>0$ such that $\nu_{a,b} (B({\bf z}, r)) =0$.
Then, by a similar argument as above, we have
\begin{align}
\nu_{a,b} (B({\bf z}, r)) 
&= \int_{\cW}  {\bf 1}_{B({\bf z}, r)} ({\bf L}(w))   \mu_{a,b}(dw)
\nn\\
&\le 
\int_{\cW} G_{{\bf z}, r} (w)   \mu_{a,b}(dw)
= \E [G_{{\bf z}, r} \,  \delta_b (Y(1,a))] /Z. 
\label{eq.zerozero}
\end{align}
Since $\bar{A}$ is closed, we can find $r>0$ such that 
$B({\bf z}, 2r) \cap \bar{A} =\emptyset$.
By \eqref{eq.indG}, $G_{{\bf z}, r}$ vanishes if ${\bf L}(w) \notin B({\bf z}, 2r)$.
Hence, we cannot find $h \in \cH$ that 
satisfies the second one of the two 
equivalent conditions in Theorem \ref{aks>0},
which implies that the right hand side of \eqref{eq.zerozero} vanishes.
\end{proof}


Under the conditions of Theorem \ref{or.supp.GRP},
we study the law of the process  $\tilde{Y}(\,\cdot\,, a)$, 
 an $\infty$-quasi-continuous modification of $Y (\,\cdot\,, a)$, 
under the probability measure $\mu_{a,b}$.
Heuristically, it is the law of ``$Y (\,\cdot\,, a)$ conditioned  
 $Y (1, a)=b$."
 Since $\tilde{Y} (\,\cdot\,, a) = \Pi \Phi^a ({\bf L}(w))$,
 the above law equals the law of $\Pi \circ \Phi^a$ under $\nu_{a,b}$.
  
 Let us make sure that this law actually sits on 
 $C_{\Pi (a), b}^{0, \beta \textrm{-H}}(\cV)$ for every $1/3 <\beta <1/2$.
 Choose $\alpha$ and $m$ so that $\beta \le \alpha- (4m)^{-1}$.
It is sufficient to show that the end point is almost surely $b$.
For every $g \in C_0^\infty (\cV)$, we have 
\begin{align*}
\int_{W}  g (\tilde{Y} (1, a)) \mu_{a,b}(dw) 
  &= \E [ g (Y (1, a)) \,  \delta_b (Y(1,a))] / Z
   \\
    &= \lim_{n\to\infty} 
      \E [ g (Y (1, a)) \,  \psi_n (Y(1,a))] / Z
      \\
    &= \lim_{n\to\infty}  \int_{\cV}  g(y) \psi_n (y) q(y) dy  /Z
      \\
    &= g(b) q(b)/Z =g(b).
           \end{align*}
   This implies that the law of $\tilde{Y} (1, a)$ under $\mu_{a,b}$
          is the point mass at $b$.
  Here, {\rm (i)} $Z:=\E [ \delta_b (Y(1,a))]>0$ is the normalizing constant,
  {\rm (ii)} $q(y)$ is the smooth density 
     (with respect to the Lebesgue measure $dy$ on $\cV$) of  
   of the law of $Y(1,a)$ under the Wiener measure, 
     and {\rm (iii)} $\{\psi_n \}_{n=1}^\infty \subset C_0^\infty (\cV)$ is any sequence 
      that converges to $\delta_b$ in the space of Schwartz distributions 
       on $\cV$.
        Note that $Z = q(b)$ due to Item {\bf (c)} in Subsection 2.1.


By a similar argument as above, we can see that 
the finite dimensional distribution of this law is uniquely
determined by the following formula:  
For every $k\in \N$,  $\{0 <t_1<\cdots < t_k <1\}$
and $g_1, \ldots, g_k \in C_0^\infty (\cV)$, it holds that
\begin{align}\label{eq.KC1}
\lefteqn{
\int_{W}  \prod_{i=1}^k g_i (\tilde{Y} (t_i, a)) \mu_{a,b}(dw) 
}
\nn\\
&=
\int_{W}  \prod_{i=1}^k g_i (\Pi \tilde{X} (t_i, a)) \mu_{a,b}(dw) 
\nn\\
  &=
 Z^{-1}
\E [\prod_{i=1}^k g_i (\Pi X (t_i, a)) \,  \delta_b (Y(1,a))]
\nn\\ 
  &=
 Z^{-1}
 \lim_{n\to\infty}\E [\prod_{i=1}^k g_i (\Pi X (t_i, a)) \,  \psi_n (Y(1,a))]
\nn\\
  &=
 Z^{-1} \lim_{n\to\infty}
  \bigl[
  T_{t_1} (g_1 \circ \Pi) T_{t_2-t_1} \cdots (g_{k-1} \circ \Pi)
   T_{t_k-t_{k-1}} (g_k \circ \Pi) 
  T_{1-t_k}  (\psi_n\circ \Pi )\bigr](a).
    \end{align}
Here, $\{\psi_n \}$ is the same as above.
Note that $g_i\circ \Pi~(1\le i \le k)$ on the right hand side 
of \eqref{eq.KC1} are viewed as 
multiplication operators.
Note also that the limit above exists and is independent of the choice of $\{\psi_n \}$.


\begin{corollary} \label{cor.supp_Y}
Let the notation and the situation be as above.
Assume \eqref{eq.non_emp} and non-degeneracy of $Y(1,a)$.
Then, the support of the law of the process  $\tilde{Y}(\,\cdot\,, a)$
under $\mu_{a,b}$ equals the closure of 
\begin{equation}\label{eq.PQR}
\{ \Pi \Psi^a (h)  \colon   \mbox{$h \in \cH$,  ~$D\Pi \Psi^a(h)_1
\colon \cH\rightarrow \mathcal{V}$ is surjective, 
~$\Pi \Psi^a(h)_1=b$}
\}
\end{equation}
in $\cC_{\Pi (a), b}^{0, \beta \textrm{-H}}(\cV)$ for every $1/3 <\beta <1/2$.
\end{corollary}

\begin{proof} 
The set in \eqref{eq.PQR} is denoted by $B$.
The set in \eqref{eq.KLM} is denoted by $A$ again.
Note that $\Pi \Psi^a(h)= \Pi \Phi^a (\cL (h))$ for every $h\in \cH$.
If $\beta \le \alpha- (4m)^{-1}$, 
then $\Pi \Phi^a =\Pi \circ \Phi^a
\colon G\Omega^{\textrm{B}}_{\alpha, 4m} ( {\mathbb R}^d)
\to \cC_{\Pi (a), b}^{0, \beta \textrm{-H}}(\cV)$ is continuous.
Thanks to this continuity, the proof is quite simple and straightforward.

If $\Pi \Psi^a (h)\in B$,
then its inverse image by $\Pi \Phi^a$ clearly intersects with $A$. 
By the continuity, the inverse image of every open neighborhood of 
$\Pi \Psi^a (h) \in B$ is an open subset of $G\Omega^{\textrm{B}}_{\alpha, 4m} ( {\mathbb R}^d)$ that intersects with $A$ 
and therefore its weight is strictly positive by Proposition \ref{or.supp.GRP}. 
This implies that $B$ is included in the support. 
So is $\bar{B}$ since the support is closed by definition.
 
Finally, we show that  this inclusion cannot be strict by 
showing $\bar{B}^c=(\bar{B})^c$ is of measure zero.
It is clear $(\Pi \Phi^a)^{-1} (\bar{B})$ and 
$(\Pi \Phi^a)^{-1} (\bar{B}^c)$ do not intersect. 
By the continuity, the former is closed, while the latter is open. 
Since $(\Pi \Phi^a)^{-1} (\bar{B}) \supset A$, the support of $\nu_{a,b}$ is
included by $(\Pi \Phi^a)^{-1} (\bar{B})$ by Theorem \ref{or.supp.GRP}. 
Hence, $\nu_{a,b} ((\Pi \Phi^a)^{-1} (\bar{B}) )=1$. This completes the proof. 
\end{proof}

Now we consider the special case $\cV =\R^e$ 
(therefore $\Pi$ is the identity map) and 
$X(t,a)$ is non-degenerate in the sense of Malliavin for 
every $a \in \R^e$ and $t \in (0,1]$.
Then, the law of $X(t,a)$ has a
 density with respect to the Lebesgue measure $db$,
 which is denoted by $p(t, a, b)$,
that is, $\mu ( X(t,a) \in A) = \int_A p(t, a, b) db$ for 
every Borel subset  $A \subset \R^e$.
In this case, the law of the process in Corollary \ref{cor.supp_Y}
is identical to the classical pinned diffusion measure ${\mathbb Q}_{a,b}$
associated with SDE \eqref{sde_x.def}. 
Indeed, the right hand side of \eqref{eq.KC1} reads:
\[
p (1, a,b)^{-1} 
\int_{(\R^e)^k}  \Bigl\{ \prod_{i=1}^k g_i (b_i) \Bigr\} p (t_1, a,b_1)  
\Bigl\{\prod_{i=2}^{k}  p (t_i - t_{i-1}, b_{i-1}, b_{i})  \Bigr\}
 p (1- t_k, b_{k},b) 
\Bigl\{\prod_{i=1}^k db_i \Bigr\}.
\]
This is the finite dimensional distribution of 
the classical pinned diffusion measure from $a\in \R^e$ to $b \in \R^e$.
Note that our argument automatically shows the existence of ${\mathbb Q}_{a,b}$. 
As a special case of the above corollary, we then have the following:

\begin{corollary} \label{cor.supp_X}
Let the notation and the situation be as above.
Assume \eqref{eq.non_emp} (with $\Pi$ being the identity map of $\R^e$)
and non-degeneracy of $X(t, x)$ for all $x \in \R^e$ and $t \in (0,1]$.
Then, the support of ${\mathbb Q}_{a,b}$ equals the closure of 
\begin{equation}\nn
\{ \Psi^a (h)  \colon   \mbox{$h \in \cH$, ~$D \Psi^a(h)_1
\colon \cH\rightarrow \mathbb{R}^e$ is surjective, 
~$\Psi^a(h)_1=b$}
\}
\end{equation}
 in $\cC_{a, b}^{0, \beta\textrm{-H}}(\R^e)$ for every $1/3 <\beta <1/2$.
\end{corollary}

\begin{remark}\label{rem.nondeg}
In this remark
we provide two typical sufficient conditions for non-degeneracy of $Y(t, a)$.  
Both are bracket-generating conditions of H\"ormander-type.

Let $V_i~(0\le i \le d)$ be the coefficient vector fields of 
SDE \eqref{sde_x.def}.
In this remark they are viewed as  first-order differential operators on $\R^e$.
Set $\Sigma_1 =\{V_1, \ldots, V_d\}$ and, recursively,
$\Sigma_k =\{ [Z, V_i] \colon  Z\in \Sigma_{k-1},  0\le i \le d\}$
for $k\ge 2$. 
\begin{enumerate}
\item[(A)] If $\{ Z(a) \colon Z \in \cup_{k=1}^\infty\Sigma_k \}$ linearly spans $\R^e$
at the starting point $a$, then for all $t \in (0,1]$ $X(t,a)$ is
non-degenerate and therefore so is $Y(t,a)$.
This fact is well-known. (See \cite[Section 2.3]{nu} or \cite[Section V.10]{iwbk}
for example.)
\item[(B)]
Suppose the following uniform partial H\"ormander condition:  
There exists $L >0$ such that 
\begin{equation}\nn
\inf_{ a \in  {\mathbb R}^e }  \inf_{\eta \in {\cal V}, |\eta| =1} 
\sum_{k=1}^L \sum_{Z \in \Sigma_k}
\la \Pi  Z (a), \eta \ra^2 >0.
\end{equation}
Then, according to \cite[Theorem 2.17 and Lemma 5.1]{ks2},
$Y(t,a)$ is non-degenerate in 
the sense of Malliavin calculus for every $a\in\R^e$ and $t\in (0,1]$.\end{enumerate}
\end{remark}

\begin{example}
We provide some examples of the process $Y(\cdot,a)$ in Corollary \ref{cor.supp_Y} (except that in Corollary \ref{cor.supp_X}).
\begin{itemize}
\item
Assume Condition (A) in Remark \ref{rem.nondeg} and $\cV=\R^e$.
Then, the solution $X(\cdot, a)$ of SDE \eqref{sde_x.def}
satisfies Corollary \ref{cor.supp_Y}.
In this case, the density $p(t, z, z')$ may not exist
if $z\in\R^e$ is  distant from $a$.
Therefore, it is not clear whether the pinned diffusion measure 
in the usual sense exist or not.
Since our method is based on quasi-sure analysis, 
we can deal with this kind of situation (without any additional efforts), too.

\item
For $1\le e' <e$, set 
$\cV= \R^{e'} \cong \R^{e'}\oplus \{{\bf 0}_{e-e'}\} \subset \R^{e}$.
Here, ${\bf 0}_{e-e'}$ is the zero vector of $\R^{e-e'}$.
If we write $X_t =(X_t^1, \ldots, X_t^e)$,
then $Y_t =\Pi X_t =(X_t^1, \ldots, X_t^{e'})$.
This kind of projected process are sometimes studied.
For example, in \cite{dfjv1, dfjv2, tawa}, small noise problems 
for the density of $Y_t$ are studied.
Therefore, it looks natural to study the pinned process conditioned by
$Y_1=b$.
(The Markov property is lost after the projection in general. 
So, it cannot be called a ``pinned diffusion process.")

\item
Assume that $V_i (t, x)\colon [0,1]\times \R^e \to \R^e$ 
extend to $C_b^\infty$-maps on an open neighborhood of 
$[0,1]\times \R^e \subset \R^{e+1}=\{(t,x)\colon t \in \R, x\in\R^e\}$  ($0\le i \le d$).
We extend them to $C_b^\infty$-maps on $\R^{e+1}$
(which will be denoted by the same symbols)
and view them as time-dependent vector fields on $\R^e$. 
Instead of \eqref{sde_x.def}, we now consider the 
following time-dependent SDE:
\[
d\tilde{X}_t = \sum_{i=1}^d  V_i (t, \tilde{X}_t)\circ dw_t^i + V_0 (t, \tilde{X}_t) dt,
\qquad
 \tilde{X}_0 =a \in \R^e.
\]
Define $\Sigma_{k}(t)$ in the same way as in Remark \ref{rem.nondeg}
by just replacing $V_i~(0\le i \le d)$ by $V_i(t, \cdot)~(0\le i \le d)$.
Some examples of bracket-generating condition 
sufficient for the non-degeneracy of $\tilde{X}_t =\tilde{X}(t,a)$
can be found in \cite{flor, tani} among others.
If we set $X_t =(X^0_t, \tilde{X}_t)$ with $X^0_t \equiv t$, 
then $X$ satisfies the following SDE on $\R^{e+1}$:
\begin{align*}  
dX_t = \sum_{i=1}^d  \hat{V}_i (X_t)\circ dw_t^i + \hat{V}_0 (X_t) dt,
\qquad
X_0 =(0,a) \in \R^{e+1}.
 \end{align*}
Here, we set $\hat{V}_0 :=V_0 + (\partial/\partial t)$,
$\hat{V}_i :=V_i$ for $1\le i \le d$ and
view them as vector fields on $\R^{e+1}$.
Since $\Pi (X_t) =\tilde{X}_t$ for the canonical projection $\Pi\colon \R^{e+1}\to \R^e$, the process $\tilde{X} (\cdot, a)$ satisfies the assumptions of Corollary \ref{cor.supp_Y}.
\end{itemize}
\end{example}

As one can easily see, 
our support theorem for pinned cases looks clearly different 
from the standard version of the support theorem
because of the two conditions on the skeleton ODE.
The first one, which is quite easy for everyone to guess, 
 requires the solution of  the skeleton ODE to end at the given point. 
The second one requires the tangent map of the solution map of 
 the skeleton ODE (at time $1$) to be non-degenerate.
This may look a little bit surprising to some readers, 
but is actually quite natural
from the viewpoint of positivity theorems for the densities for SDEs.

%

\medskip
\noindent
{\bf Acknowledgement:}~
The author is supported by 
JSPS KAKENHI (Grant No. 20H01807).


\bigskip
\begin{flushleft}
  \begin{tabular}{ll}
    Yuzuru \textsc{Inahama}
    \\
    Faculty of Mathematics,
    \\
    Kyushu University,
    \\
    744 Motooka, Nishi-ku, Fukuoka, 819-0395, JAPAN.
    \\
    Email: {\tt inahama@math.kyushu-u.ac.jp}
  \end{tabular}
\end{flushleft}

\end{document}